\input amstex
\input diagrams
\documentstyle{amsppt}
\document
\topmatter
\title
On the Eisenstein ideal of Drinfeld modular curves
\endtitle
\author
Ambrus P\'al
\endauthor
\date
April 14, 2006.
\enddate
\abstract Let $\goth E(\goth p)$ denote the Eisenstein ideal in the Hecke algebra $\Bbb T(\goth p)$ of the Drinfeld modular curve $X_0(\goth p)$ parameterizing Drinfeld modules of rank two over $\Bbb F_q[T]$ of general characteristic with Hecke level $\goth p$-structure, where $\goth p\triangleleft\Bbb F_q[T]$ is a non-zero prime ideal. We prove that the characteristic $p$ of the field $\Bbb F_q$ does not divide the order of the quotient $\Bbb T(\goth p)/\goth E(\goth p)$ and the Eisenstein ideal $\goth E(\goth p)$ is locally principal.
\endabstract
\footnote" "{\it 2000 Mathematics Subject Classification. \rm 11G18
(primary), 11G09 (secondary).}
\endtopmatter

\heading 1. Introduction
\endheading

\definition{Notation 1.1} Let $F=\Bbb F_q(T)$ denote the rational function
field of transcendence degree one over a finite field $\Bbb F_q$ of
characteristic $p$, where $T$ is an indeterminate, and let $\bold A=\Bbb
F_q[T]$. Let $Y_0(\goth n)$ be the Drinfeld modular curve parameterizing
Drinfeld modules of rank two over $\Bbb F_q[T]$ of general characteristic with Hecke level $\goth n$-structure, where $\goth n\triangleleft\Bbb F_q[T]$ is a non-zero ideal. Let $X_0(\goth n)$ denote the unique geometrically irreducible non-singular projective curve containing $Y_0(\goth n)$ and let moreover $J_0(\goth n)$ denote the Jacobian of the curve $X_0(\goth n)$. Assume now that $\goth n=\goth p$ is a prime ideal and let $\goth E(\goth p)$ denote the Eisenstein ideal
in the Hecke algebra $\Bbb T(\goth p)$ of the Drinfeld modular curve $X_0(\goth p)$, defined in [22] first in this context. The latter was already studied intensively in [19] already. There a close analogue of the classical theory in [16] was worked out in detail. In this paper we will complete our previous results. In particular we show:
\enddefinition
\proclaim{Theorem 1.2} We have:
$${\Bbb T(\goth p)\over\goth E(\goth p)}={\Bbb Z\over N(\goth p)\Bbb Z},$$
where
$$N(\goth p)=\cases{q^{\deg(\goth p)}-1\over q-1},&\text{if $\deg(\goth p)$ is odd,}\\{q^{\deg(\goth p)}-1\over q^2-1},&\text{otherwise.}\endcases$$
\endproclaim
By the results of [19] the new information in the claim above is that the characteristic $p$ does not divide the order $|\Bbb T(\goth p)/\goth E(\goth p)|$. The latter has the following immediate
\proclaim{Corollary 1.3} There is no non-zero double-cuspidal, weight two Drinfeld modular form $\phi\in M^2_{2,1}(\Gamma_0(\goth p))$ of type $1$ and
level $\goth p$ fixed by the Hecke algebra $\Bbb T(\goth p)$.
\endproclaim
In fact we will prove a similar claim about $M^2_{2,1}(\Gamma_0(\goth n))$ for every square-free ideal $\goth n\triangleleft\bold A$. The other main result of this paper is the following
\proclaim{Theorem 1.4} The Eisenstein ideal $\goth E(\goth p)$ is locally principal.
\endproclaim
We are going to give two proofs of the result above. In the first proof theorem is
deduced from the fact that the ring $\Bbb T(\goth p)$ is locally Gorenstein at the prime ideals lying above the Eisenstein ideal, which we proved already in [19], by adopting Mazur's original argument to the modular symbols introduced in [23]. In fact we will produce explicit generators (for details see Theorem 4.17). In the second proof we will identify (a suitable enlargement of) $\Bbb T(\goth p)$ with a universal deformation ring $R(\goth p)$ analogous to the construction in [1] using the Wiles-Lenstra criterion, then prove that $R(\goth p)$ is generated by one element over $\Bbb Z_l$ using cohomological methods. The deformation theoretical methods prove directly the Gorenstein property without using much about the finer geometry of the modular curve $X_0(\goth p)$ hence also give an alternative route to prove the main diophantine results of the paper [19].
\definition{Contents 1.5} In the next chapter we are going to prove that every abelian variety defined over $F$ whose N\'eron model has either good ordinary reduction or multiplicative reduction at every closed point of the projective line over $\Bbb F_q$ has everywhere good reduction. In the third chapter we prove some preliminary results about harmonic cochains. The result proved in the second chapter is used in the fourth chapter to prove Theorem 1.2 and Corollary 1.3 where we also define the Eisenstein ideal of the Jacobian of the Drinfeld modular curve $X_0(\goth n)$. In the fifth chapter we develop theory of modular symbols for Drinfeld modular curves, first studied by Teitelbaum in this context, and use it to give the first proof of Theorem 1.4. We introduce a deformation functor analogous to the one defined in [1] and prove that it is representable in the sixth chapter. In the seventh chapter we construct a homomorphism from the universal deformation ring of the functor introduced in the previous chapter into a certain Hecke algebra. We prove that this map is an isomorphism in the eighth chapter and deduce the second proof of Theorem 1.4 from this result. The aim of the ninth chapter is to deduce the main diophantine results of [19] directly from the Gorenstein property.
\enddefinition
\definition{Acknowledgment 1.6} I wish to thank the IH\'ES for
its warm hospitality and the pleasant environment they created for
productive research, where this article was written.
\enddefinition

\heading 2. Abelian schemes with multiplicative or ordinary reduction everywhere
\endheading

\definition{Definition 2.1} Let $\pi:A\rightarrow X$ be an abelian scheme of relative dimension $d$ defined over an $\Bbb F_p$-scheme $X$, where $p$ is a prime number. Let $\pi^{(p^n)}:A^{(p^n)}\rightarrow X$ be the pull-back of $A$ with respect to the $n$-iterated absolute Frobenius $F^n:X\rightarrow X$ and let $F^n:A\rightarrow A^{(p^n)}$ be the $n$-iterated relative Frobenius. Then $A^{(p^n)}$ is an abelian scheme over $X$ and $F^n$ is a morphism of group schemes. We say that $A$ is ordinary if the $\bold k$-valued points of $A[p]$ form a group of order $p^d$ for every geometric point $\bold k$ of $X$. A finite flat group scheme $G$ over $X$ is called multiplicative if its Cartier dual is an \'etale group scheme.
\enddefinition
\proclaim{Proposition 2.2} Let $A$ be an ordinary abelian scheme of dimension $d$ over a noetherian $\Bbb F_p$-scheme $X$. Then $F^n$ is flat, surjective and finite, and its kernel $Ker(F^n)$ is a multiplicative finite flat group scheme of rank $p^{nd}$ over $X$.
\endproclaim
\definition{Proof} This proposition is certainly well-known at least when $X$ is the spectrum of a field. We deduce the general case from the above in the usual manner. The map $F^n$ is flat because $A$ is flat over $X$ and $F^n$ is flat on each fiber of $A$ over $X$ by Theorem 5.9 of section IV of [13], pages 99-100, known as the local criterion for flatness. The composition $\pi=\pi^{(p^n)}\circ F^n$ and the map $\pi^{(p^n)}$ are both proper and separated, hence $F^n$ is proper, too. This map has finite fibers, so it is a finite map by Corollary 1.10 of [18], pages 6-7. As $Ker(F^n)$ is the fiber of $F^n$ over the zero section we get that it is finite and flat over $X$. On the other hand it is also a group scheme by definition. In order to prove the two remaining claims about the order and the multiplicativity of $Ker(F^n)$, it will be sufficient to show that the $\bold k$-valued points of the Cartier dual of $Ker(F^n)$ form a group of order $p^{nd}$ for every geometric point $\bold k$ of $X$ as the formation of the Cartier dual commutes with base change. The latter is just a reformulation of the claims mentioned above when $X$ is the spectrum of an algebraically closed field.\ $\square$
\enddefinition
\definition{Definition 2.3} Let $C$ denote a smooth irreducible curve defined over a field $\bold k$ of characteristic $p$. For every closed point $x$ of $C$ let $\Cal O_x$ denote the local ring of $C$ at $x$. Moreover let $\widehat{\Cal O}_x$, $K_x$ denote the completion of $\Cal O_x$ and the quotient field of $\widehat{\Cal O}_x$, respectively. Let $K$ be the function field of $C$ and let $A$ be an abelian variety defined over $K$. Finally let $\Cal A$ denote the N\'eron model of $A$ over $C$. Recall that $A$ has multiplicative reduction at a closed point $x$ if the connected component of the identity of the fiber of $\Cal A$ over $x$ is a torus.
\enddefinition
\proclaim{Lemma 2.4} Assume that $A$ has multiplicative reduction at the closed point $x$ of $C$. Then $Ker(F^n)$ in $A$ extends to a multiplicative finite flat group scheme over $\Cal O_x$ of rank $p^{nd}$.
\endproclaim
\definition{Proof} It will be enough to show that the Cartier dual of $Ker(F^n)$ extends to a finite \'etale group scheme over $\Cal O_x$. In other words we need to show that the Galois representation attached to this \'etale group scheme is unramified. Since the formation of the Cartier dual commutes with base change, it will be enough to show that $Ker(F^n)$ as a group scheme over Spec$(K_x)$ extends to a multiplicative finite flat group scheme over $\widehat{\Cal O}_x$ of the prescribed rank. Since $A$ has multiplicative reduction by assumption, it is semi-stable, and its Raynaud group, which is a smooth group scheme $\Cal A^{\#}$ over $\widehat{\Cal O}_x$ of finite type in general, is a torus. The formal completions of $\Cal A^{\#}$ and $\Cal A$ along the special fiber are equal. The finite flat group scheme $\Cal A^{\#}[p^n]$ is multiplicative of rank $p^{nd}$ and annihilated by $F^n$, hence $\Cal A[p^n]$ also contains such a group scheme. The base change of the latter to Spec$(K_x)$ lies in $Ker(F^n)$, but their rank is the same, so they must be equal.\ $\square$
\enddefinition
\definition{Definition 2.5} We will continue to use the notation above. For any category $\Cal C$ let $\text{Ob}(\Cal C)$ denote the class of objects of $\Cal C$. For any scheme $X$ let $\bold{FFGp}(X)$ denote the category of finite flat group schemes over $X$. For any morphism of schemes $m:Y\rightarrow X$ let $m^*:\bold{FFGp}(X)\rightarrow \bold{FFGp}(Y)$ denote the functor induced by the pull-back with respect to $m$. Let $S$ be a finite set of closed points of $C$ and let $X$ denote the complement of $S$ in $C$. Let $j:\text{Spec}(F)\rightarrow X$ denote the generic point of $X$ and let $i:X\rightarrow C$ be its closed immersion into $C$. For every $x$ closed point of $C$ let $j_x:\text{Spec}(F)\rightarrow\text{Spec}(\Cal O_x)$ and $i_x:\text{Spec}(\Cal O_x)\rightarrow C$ be the generic point and closed immersion of the corresponding scheme, respectively.  Let $\bold{DD}(X,S)$ denote the additive category whose objects are collections of finite flat group schemes $G_X\in\text{Ob}(\bold{FFGp}(X))$ and $G_x\in\text{Ob}(\bold{FFGp}(\text{Spec}(\Cal O_x))$ for every $x\in S$ and an isomorphism $g_x:j^*(G_X)\rightarrow j_x^*(G_x)$ for every $x\in S$, and a morphism $\phi:(G_X,G_x,g_x)\rightarrow(H_X,H_x,h_x)$ in this category is a collection of morphisms $\phi_X:G_X\rightarrow H_X$ and $\phi_x:G_x\rightarrow H_x$ for every $x\in S$ such that $j_x^*(\phi_x)\circ g_x=h_x\circ j^*(\phi_X)$. Since $j_x\circ i_x=j\circ i$ for every $x$ closed point of $C$, there is a natural isomorphism $c_x$ between the functors $j_x^*\circ i_x^*$ and $j^*\circ i^*$. Let $d$ denote the functor $d:\bold{FFGp}(C)\rightarrow\bold{DD}(X,S)$ which assigns the object $G\in\text{Ob}(\bold{FFGp}(C))$ the collection $(i^*(G),i_x^*(G),c_x(G))$ and assigns the morphism $\phi:G\rightarrow H$ the collection $(i^*(\phi),i_x^*(\phi))$. This is clearly a faithful embedding of categories.
\enddefinition
\proclaim{Lemma 2.6} The functor $d$ is an equivalence of categories.
\endproclaim
\definition{Proof} The proof is an exercise in elementary descent theory. First we need to show that every object $(G_X,G_x,g_x)$ of $\bold{DD}(X,S)$ is in fact isomorphic to an object of $\bold{FFGp}(C)$. Since every finite and flat morphism is in fact affine, there is a coherent sheaf of Hopf algebras $A_X$ over $X$ and a coherent sheaf of Hopf algebras $A_x$ over a Zariski neighborhood of $x$ for every $x\in S$ whose spectrum is $G_X$ over $X$ and $G_x$ over Spec$(\Cal O_x)$, respectively. These sheaves patch together to a coherent sheaf of Hopf algebras over $C$ via the maps induced by the morphisms $g_x$ whose spectrum is the finite flat group scheme $G$ we were looking for. A similar argument shows that any morphism in $\bold{DD}(X,S)$ comes from a morphism in $\bold{FFGp}(C)$.\ $\square$
\enddefinition
Now we assume that $\bold k$ is a finite field and $C$ is the projective line $\Bbb P^1_{\bold k}$ over $\bold k$.
\proclaim{Theorem 2.7} Let $A$ be an abelian variety defined over $K$ such that the N\'eron model $\Cal A$ of $A$ has either good ordinary reduction or multiplicative reduction at every closed point of $C$. Then $A$ has everywhere good reduction.
\endproclaim
\definition{Proof} We may assume that there is a degree one closed point $c$ in $C$ where $A$ has good reduction by enlarging $\bold k$ if necessary. Let $\Cal B$ denote the unique constant abelian group scheme over $C$ whose fiber $\Cal B_c$ at $c$ is isomorphic to the fiber $\Cal A_c$ of $\Cal A$ at $c$. Let $S$ denote the set of closed points of $C$ where $A$ has multiplicative reduction and let $X$ denote the complement of $S$ in $C$ as above. By Proposition 2.2 the kernel of $F^n$ in the abelian scheme $\Cal A|_X$ is a multiplicative finite flat group scheme of rank $p^{nd}$ over $X$ where $d$ is the dimension of $A$ over $K$. On the other hand the kernel of $F^n$ in the abelian variety $A$ extends to a multiplicative finite flat group scheme over $\Cal O_x$ of rank $p^{nd}$ for every $x\in S$ by Lemma 2.4. These group schemes glue together to a finite flat group scheme $A_n$ over $C$ by Lemma 2.6. Since the property of being multiplicative is local with respect to the Zariski topology, this group scheme is multiplicative. Therefore its Cartier dual is \'etale. Since the geometric fundamental group of $C$ is trivial, we get that $A_n$ is the base change of a group scheme over Spec$(\bold k)$ with respect to the constant map $C\rightarrow\text{Spec}(\bold k)$. Its fiber over $c$ is isomorphic to the fiber of the kernel of $F^n$ in $\Cal B|_X$ at $c$, so these group schemes must be isomorphic. Hence  the fibers $\Cal A_x$ and $\Cal B_x$ must be isogenous for every closed point $x$ in $X$ by Lemma 2.8 below. Therefore $A$ and $B$ must be isogenous by Zarhin's isogeny theorem where $B$ is the generic fiber of $\Cal B$. In particular $A$ has everywhere good reduction.\ $\square$
\enddefinition
\proclaim{Lemma 2.8} Let $\bold f$ be a finite field of characteristic $p$ and let $C_1$ and $C_2$ be two ordinary abelian varieties of same dimension defined over $\bold f$. Assume that the kernel of $F^n$ in $C_1$ and in $C_2$ are isomorphic as group schemes for every positive integer $n$. Then $C_1$ and $C_2$ are isogenous over $\bold f$.
\endproclaim
\definition{Proof} Let $\Bbb W(\bold f)$ and $\Bbb Q(\bold f)$ denote the Witt vectors of infinite length over $\bold f$ and its quotient field, respectively. Let $\sigma:\Bbb Q(\bold f)\rightarrow \Bbb Q(\bold f)$ denote the Frobenius. For any finite flat group scheme $G$ over $\bold f$ let $\Bbb D(G)$ denote its contravariant Dieudonn\'e module. Because $C_1$ and $C_2$ are ordinary, there are exact sequences:
$$0@>>>M_{i,n}@>>>C_i[p^n]@>>>N_{i,n}@>>>0,$$
where $M_{i,n}$, $N_{i,n}$ are multiplicative and \'etale group schemes of rank $p^{nd}$, respectively, for $i=1$ or $i=2$ and every $n\in\Bbb N$. Recall that the Dieudonn\'e crystal $\Bbb D(V)$ associated to an abelian variety $V$ over $\bold f$ is defined as the limit of the Dieudonn\'e modules $\Bbb D(V[p^n])$ with respect to the maps $\Bbb D(V[p^m])\rightarrow\Bbb D(V[p^n])$ induced by the closed immersions $V[p^n]\rightarrow V[p^m]$ for every $n\leq m$. By our assumption the group schemes $M_{1,n}$ and $M_{2,n}$ are isomorphic for every $n\in\Bbb N$. Because there are only finitely many isomophisms between $M_{1,n}$ and $M_{2,n}$, we may chose one $\phi_n:M_{1,n}\rightarrow M_{2,n}$ for every $n$ such that $\phi_n|_{M_{1,m}}=\phi_m$ for every $m\leq n$ by Tychonov's theorem. Hence the exact sequence above induces an exact sequence:
$$0@>>>\Bbb N_i@>>>\Bbb D(C_i)@>>>\Bbb M@>>>0$$
for $i=1$ or $i=2$, where $\Bbb M$ is a Dieudonn\'e crystal of slope $1$ and $\Bbb N_i$ are Dieudonn\'e crystals of slope $0$. Recall that for every $F$-crystal $\Bbb V=(V,\phi)$ over $\Bbb Q(\bold f)$, where $V$ is a finite dimensional vector space over $\Bbb Q(\bold f)$ and $\phi:V\rightarrow V$ is a bijective $\sigma$-linear map, its dual $\Bbb V^{\vee}$ is defined as the pair $(\text{Hom}_{\Bbb Q(\bold f)}
(V,\Bbb Q(\bold f)),\phi^{\vee})$, where $\phi^{\vee}(f)(u)=pf(\phi^{-1}(u))^{\sigma}$ for every $f\in\text{Hom}_{\Bbb Q(\bold f)}(V,\Bbb Q(\bold f))$ and $u\in V$. Every choice of polarization of $C_i$ induces an isomorphism of $F$-crystals $(\Bbb D(C_i)\otimes_{\Bbb W(\bold f)}\Bbb Q(\bold f))^{\vee}=\Bbb D(C_i)\otimes_{\Bbb W(\bold f)}\Bbb Q(\bold f)$ over $\Bbb Q(\bold f)$. Such an isomorphism induces an isomorphism between the largest  quotient modules of slope $1$ which are equal to $(\Bbb M\otimes_{\Bbb W(\bold f)}\Bbb Q(\bold f))^{\vee}$ and $\Bbb N_i\otimes_{\Bbb W(\bold f)}\Bbb Q(\bold f)$, respectively. Therefore the semi-simplifications of $\Bbb D(C_i)\otimes_{\Bbb W(\bold f)}\Bbb Q(\bold f)$, where $i=1$ or $i=2$, are isomorphic. Hence the characteristic polynomial of the arithmetic Frobenius of the Dieudonn\'e crystals $\Bbb D(C_1)$ and $\Bbb D(C_2)$ are equal which implies that $C_1$ and $C_2$ must be isogenous by Honda-Tate theory.\ $\square$
\enddefinition

\heading 3. Harmonic cochains
\endheading

\definition{Definition 3.1} For any graph $G$ let $\Cal V(G)$ and $\Cal
E(G)$ denote its set of vertices and edges, respectively. In this paper
we will only consider such oriented graphs $G$ which are equipped with an involution $\overline{\cdot}:\Cal E(G)\rightarrow\Cal E(G)$ such that for each edge $e\in\Cal E(G)$ the original and terminal vertices of the edge $\overline e\in\Cal E (G)$ are the terminal and original vertices of $e$, respectively. The edge $\overline e$ is called the edge $e$ with reversed orientation. Let $R$ be a commutative group. A function $\phi:\Cal E(G)\rightarrow R$ is called a harmonic $R$-valued cochain, if it satisfies the following conditions:
\roster
\item"$(i)$" We have:
$$\phi(e)+\phi(\overline e)=0\text{\ }(\forall e\in\Cal E(G)).$$
\item"$(ii)$" If for an edge $e$ we introduce the notation $o(e)$ and
$t(e)$ for its original and terminal vertex respectively,
$$\sum_{{e\in\Cal E(G)\atop o(e)=v}}\phi(e)=0\text{\ }
(\forall v\in\Cal V(G)).$$
\endroster
We denote by $H(G,R)$ the group of $R$-valued harmonic cochains on $G$.
\enddefinition
\definition{Definition 3.2} Let $GL_2$ denote the group scheme of invertible two by two matrices and let $Z$ denote its center. Let $F_{\infty}$ denote the completion of $F=\Bbb F_q(T)$ with respect to the valuation $\infty$ corresponding to the point at infinity on the projective line over $\Bbb F_q$. Let $\Cal O_{\infty}$ denote the valuation ring of $F_{\infty}$ and let $\upsilon\in F_{\infty}$ be a uniformizer. We are going to recall the definition of the Bruhat-Tits tree $\Cal T$ associated to the projective linear group $PGL_2(F_{\infty})$. The set of vertices $\Cal V(\Cal T)$ and edges $\Cal E(\Cal T)$ are the cosets $GL_2(F_{\infty})/GL_2(\Cal O_{\infty})Z(F_{\infty})$ and $GL_2(F_{\infty})/\Gamma_{\infty}Z(F_{\infty})$, respectively, where $\Gamma_{\infty}$ is the Iwahori group:
$$\Gamma_{\infty}=\left\{\left(\matrix a&b\\
c&d\endmatrix\right)\in GL_2(\Cal O_{\infty})|
\infty(c)>0\right\}.$$
Since $\Gamma_{\infty}$ is a subgroup of $GL_2(\Cal O_{\infty})$ there is a natural morphism $o:\Cal E(\Cal T)\rightarrow\Cal V(\Cal T)$ which assigns to every edge its original vertex. The matrix $\left(
\smallmatrix0&1\\\upsilon&0\endsmallmatrix\right)$ normalizes the Iwahori subgroup therefore the map $GL_2(F_{\infty})\rightarrow GL_2(F_{\infty})$
given by the rule $\left(\smallmatrix a&b\\c&d\endsmallmatrix\right)
\rightarrow\left(\smallmatrix a&b\\c&d\endsmallmatrix\right)
\left(\smallmatrix0&1\\\upsilon&0\endsmallmatrix\right)$
induces a map on the coset $\Cal V(\Cal T)$. This map is the involution which assigns to every edge $e$ the same edge $\overline e$ with reversed orientation. The composition of this involution and the map $o$ is the map $t:
\Cal E(\Cal T)\rightarrow\Cal V(\Cal T)$ which assigns to every edge its terminal vertex.
\enddefinition
\definition{Definition 3.3} For every non-zero ideal $\goth n\triangleleft\bold A$ let $\Gamma_0(\goth n)$ denote the Hecke congruence group:
$$\Gamma_0(\goth n)=\{\left(\matrix a&b\\c&d\endmatrix\right)
\in GL_2(\bold A)|c\equiv0\text{\ mod }\goth n\}.$$
The group $GL_2(F_{\infty})$ acts on itself via its left-regular action which induces an action of $GL_2(F_{\infty})$ on the Bruhat-Tits tree. This action induces an action of its subgroup $\Gamma_0(\goth n)$ on $\Cal T$ as well. Let $H(\Cal T,R)^{\Gamma_0(\goth n)}$ denote the group of $\Gamma_0(\goth n)$-invariant $R$-valued cochains on $\Cal T$. The group $GL_2(\bold A)$ does not contain elements which map an edge $e\in\Cal E(\Cal T)$ to the same edge $\overline e$ with reversed orientation therefore the cosets $\Gamma_0(\goth n)\backslash\Cal V(\Cal T)$ and $\Gamma_0(\goth n)\backslash\Cal E(\Cal T)$ are the vertices and edges of an oriented graph which is going to be denoted by $\Gamma_0(\goth n)\backslash\Cal T$. Every element $\phi$ of $H(\Cal T,R)^{\Gamma_0(\goth n)}$ induces an $R$-valued function on the edges of $\Gamma_0(\goth n)\backslash\Cal T$. If this function is zero outside of a finite set we say that $\phi$ is cuspidal. The $R$-module of cuspidal elements of $H(\Cal T,R)^{\Gamma_0(\goth n)}$ is denoted by $H_!(\Cal T,R)^{\Gamma_0(\goth n)}$. Finally let $H_{!!}(\Cal T,R)^{\Gamma_0(\goth n)}$ denote the image of the canonical map $H_!(\Cal T,\Bbb Z)^{\Gamma_0(\goth n)}\otimes R\rightarrow H_!(\Cal T,R)^{\Gamma_0(\goth n)}$.
\enddefinition
\definition{Definition 3.4} For every pair $\goth m$, $\goth n\triangleleft\bold A$ of non-zero ideals let $H(\goth m,\goth n)$ denote the set:
$$H(\goth m,\goth n)=
\{\left(\matrix a&b\\c&d\endmatrix\right)\in GL_2(F)|a,b,c,d\in\bold A,(ad-cb)=\goth m,\goth n\supseteq(c),(d)+\goth n=\bold A\}.$$
The set $H(\goth m,\goth n)$ is finite and it is a double $\Gamma_0
(\goth n)$-coset, so it is a disjoint union of finitely many left
$\Gamma_0(\goth n)$-cosets. Let $R(\goth m,\goth n)$ be a set of
representatives of these cosets. For any left $\Gamma_0(\goth n)$-invariant $R$-valued function $\phi:\Cal E(\Cal T)\rightarrow R$ define $T_{\goth m}(\phi)$ by the formula:
$$T_{\goth m}(\phi)(g)=\sum_{h\in R(\goth m,\goth n)}\phi(hg),\quad\forall g\in\Cal E(\Cal T).$$
It is well-known and easy to check that $T_{\goth m}(\phi)$ is independent of the choice of $R(\goth m,\goth n)$ and it is also a left $\Gamma_0(\goth n)$-invariant $R$-valued function so we have an $R$-linear operator $T_{\goth m}$ acting on the $R$-module of left $\Gamma_0(\goth n)$-invariant $R$-valued functions on $\Cal E(\Cal T)$. It is also well-known and not to difficult to verify that $T_{\goth m}$ leaves the submodules $H(\Cal T,R)^{\Gamma_0(\goth n)}$ and $H_!(\Cal T,R)^{\Gamma_0(\goth n)}$ invariant.
\enddefinition
\definition{Notation 3.5} Let $\Omega$ denote the rigid analytic
upper half plane, or Drinfeld's upper half plane over $F_{\infty}$. The
set of points of $\Omega$ is $\Bbb C_{\infty}-F_{\infty}$, denoted also by
$\Omega$ by abuse of notation, where $\Bbb C_{\infty}$ is the completion
of the algebraic closure of $F_{\infty}$. For the definition of its rigid
analytic structure as well as the other concepts recalled below see for
example [12]. For each holomorphic function $u:\Omega@>>>\Bbb
C^*_{\infty}$ let $r(u):\Cal E(\Cal T)@>>>\Bbb Z$ denote the van der
Put logarithmic derivative of $u$ (see [12], page 40). The group
$GL_2(\bold A)$ acts on Drinfeld's upper half plane $\Omega$ on the left via M\"obius
transformations. This action is discrete hence the set $\Gamma_0(\goth n)\backslash
\Omega$ is equipped naturally with the structure of a rigid analytic curve. Let $Y_0(\goth n)$ also denote the underlying rigid analytical space of the base change of
$Y_0(\goth n)$ to $F_{\infty}$ by abuse of notation.
\enddefinition
\proclaim{Theorem 3.6} There is a rigid-analytical isomorphism:
$$Y_0(\goth n)\cong \Gamma_0(\goth n)\backslash\Omega.$$
\endproclaim
\definition{Proof} This is a special case of Theorem 6.6 of [4].\ $\square$
\enddefinition
\definition{Notation 3.7} If $\psi:\bold A@>>>\Bbb
C_{\infty}\{\tau\}$ is a Drinfeld module of rank two over $\bold A$, then
$$\psi(T)=T+g(\psi)\tau+\Delta(\psi)\tau^2,$$
where $\Delta$ is the Drinfeld discriminant function. It is a Drinfeld
modular form of weight $q^2-1$. Under the identification of Theorem 3.5
the Drinfeld discriminant function $\Delta$ is a nowhere vanishing holomorphic
function on $\Omega$. For every ideal $\goth
n=(n)\triangleleft\bold A$ let $\Delta_{\goth n}$ denote the modular form of
weight $q^2-1$ given by the formula $\Delta_{\goth n}(z)=\Delta(nz)$.
As the notation indicates $\Delta_{\goth n}$ is independent of the choice
of the generator $n\in\goth n$. Let $E_{\goth n}=
r(\Delta/\Delta_{\goth n})$. Since $\Delta/\Delta_{\goth n}$ is a modular
form of weight zero, i.e. it is a modular unit, the function $E_{\goth n}$
is a $\Bbb Z$-valued harmonic cochain invariant under the action of $\Gamma_0(\goth n)$ because the van der Put derivative is equivariant with respect to the action of $GL_2(\bold A)$. For the rest of this chapter assume that $\goth n=\goth p$ is a proper prime ideal. Let $\pi\in\bold A$ be the unique monic polynomial generating $\goth p$.
For every left $\Gamma_0(\goth p)$-invariant $R$-valued function $\phi:\Cal E(\Cal T)\rightarrow R$ define $W_{\goth p}(\phi)$ by the formula:
$$W_{\goth p}(\phi)(g)=\phi(\left(\matrix0&1\\\pi&0\endmatrix\right)g),\quad
\forall g\in\Cal E(\Cal T).$$
Because the matrix in the formula above normalizes $\Gamma_0(\goth p)$ the function $W_{\goth p}(\phi)$ is also left $\Gamma_0(\goth p)$-invariant.
\enddefinition
\proclaim{Lemma 3.8} We have:
$$W_{\goth p}(E_{\goth p})=-E_{\goth p}\text{\ and\ }T_{\goth q}(E_{\goth
p}) =(1+q^{\deg(\goth q)})E_{\goth p}$$ for every prime ideal $\goth
q\triangleleft\bold A$ different from $\goth p$.
\endproclaim
\definition{Proof} This is Lemma 6.2 of [19] on page 153-154.\ $\square$
\enddefinition
\definition{Definition 3.9} Let $\mu_{\infty}$ be a Haar measure on the additive group of $F_{\infty}$ and for every $g\in GL_2(F_{\infty})$ let the same symbol denote the edge of $\Cal T$ corresponding to the coset of $g$ in $GL_2(F_{\infty})/\Gamma_{\infty}
Z(F_{\infty})$. For every left $\Gamma_0(\goth p)$-invariant complex-valued function $\phi$ on $\Cal E(\Cal T)$ the integral
$$\phi^0=\int_{\bold A\backslash F_{\infty}}\phi(\left(\matrix
1&x\\0&1\endmatrix\right))\text{d}\mu_{\infty}(x)$$
is well-defined as the integrand is invariant under translation by elements of $\bold A$ and the domain of integration is compact.
\enddefinition
\proclaim{Proposition 3.10} A harmonic cochain $\phi\in H(\Cal T,\Bbb C)^{\Gamma_0(\goth p)}$ is cuspidal if any only if the integrals $\phi^0$ and $W_{\goth p}(\phi)^0$ are both zero.
\endproclaim
\definition{Proof} This is Proposition 6.3 of [19] on page 154-155.\ $\square$
\enddefinition
We will also need the following
\proclaim{Lemma 3.11} The group $H(\Cal T,\Bbb C)^{GL_2(\bold A)}$ is trivial.
\endproclaim
\definition{Proof} Follows from Proposition 3.8 of [19] on pages 140-141.\ $\square$
\enddefinition
\proclaim{Theorem 3.12} We have:
$$H(\Cal T,\Bbb C)^{\Gamma_0(\goth p)}=H_!(\Cal T,\Bbb C)^{\Gamma_0(\goth p)}
\bigoplus\Bbb CE_{\goth p}.$$
\endproclaim
\definition{Proof} According to Proposition 5.8 of [19] on pages 151-152 the integral
$E_{\goth p}^0$ is not zero hence $E_{\goth p}$ is not an element of $H_!(\Cal T,\Bbb C)^{\Gamma_0(\goth p)}$ by Proposition 3.10 above. On the other hand $E_{\goth p}^0=-W_{\goth p}(E_{\goth p})^0$ by Lemma 3.8. Hence it will be sufficient to prove that the same identity holds for every harmonic cochain $\phi\in H(\Cal T,\Bbb C)^{\Gamma_0(\goth p)}$ by Proposition 3.10. Let $R(\goth p)\subset\Bbb F_q[T]$ denote the set of polynomials whose degree is less then $\deg(\goth p)$ and let $S(\goth p)$ denote the set
$$S(\goth p)=\{\left(\matrix1&0\\0&1\endmatrix\right)\}
\cup\{\left(\matrix0&1\\1&r\endmatrix\right)|r\in R(\goth p)\}.$$
Then $S(\goth p)$ is a set of representatives of the left $\Gamma_0(\goth p)$-cosets of $GL_2(\bold A)$ hence the function:
$$g\mapsto\sum_{s\in S(\goth p)}\phi(sg),\quad\forall g\in\Cal E(\Cal T)$$
is an element of $H(\Cal T,\Bbb C)^{GL_2(\bold A)}$. Therefore this function is zero by Lemma 3.11. In particular the integral:
$$\split\int_{\bold A\backslash F_{\infty}}\sum_{s\in S(\goth p)}
\phi(s\left(\matrix1&x\\0&1\endmatrix\right))\text{d}\mu_{\infty}(x)=&
\phi^0+\int_{\bold A\backslash F_{\infty}}
\!\sum_{r\in R(\goth p)}\!\!W_{\goth p}(\phi)
(\left(\matrix{1\over\pi}&{x+r\over\pi}\\0&1\endmatrix\right))
\text{d}\mu_{\infty}(x)\\
=&\phi^0+\int_{\pi\bold A\backslash F_{\infty}}
\!\!\!W_{\goth p}(\phi)(\left(\matrix{1\over\pi}&{x\over\pi}\\0&1\endmatrix\right))
\text{d}\mu_{\infty}(x)\endsplit$$
is zero. Since $W_{\goth p}(\phi)$ is a harmonic cochain, it
satisfies the identity:
$$W_{\goth p}(\phi)(g)=\sum_{\epsilon\in\Bbb F_q}W_{\goth p}(\phi)
(g\left(\matrix\upsilon&\epsilon\\0&1\endmatrix\right)),\quad\forall g\in GL_2(F_{\infty})$$
where $\upsilon\in F_{\infty}$ is a uniformizer. By a $\deg(\goth p)$-fold application of this identity we get the formula:
$$\split W_{\goth p}(\phi)^0=&
\int_{A\backslash F_{\infty}}\!\sum_{\epsilon\in\Bbb F_q}
W_{\goth p}(\phi)
\left(\matrix\upsilon&x+\epsilon\\0&1\endmatrix\right)\text{d}
\mu_{\infty}(x)\\=&\int_{\upsilon^{-1}A\backslash F_{\infty}}
\!\!\!W_{\goth p}(\phi)\left(\matrix
\upsilon&\upsilon x\\0&1\endmatrix\right)\text{d}\mu_{\infty}(x)\\=&\ldots=
\int_{\pi\bold A\backslash F_{\infty}}
\!\!\!W_{\goth p}(\phi)(\left(\matrix{1\over\pi}&{x\over\pi}\\0&1\endmatrix\right))
\text{d}\mu_{\infty}(x).\endsplit$$
The claim is now clear.\ $\square$
\enddefinition
\proclaim{Corollary 3.13} The $\Bbb Z$-module $H(\Cal T,\Bbb Z)^{\Gamma_0(\goth p)}$ is finitely generated and free.
\endproclaim
\definition{Proof} We are going to show that for every $\Bbb C$-subspace $W$ of finite dimension $n$ of the space of $\Bbb C$-valued functions on a set $X$ there is a subset $S\subset X$ of cardinality $n$ such that every element of $W$ is uniquely determined by its restriction to $S$. This claim clearly implies the proposition above by Theorem 3.12 as the vector space $H_!(\Cal T,\Bbb C)^{\Gamma_0(\goth p)}$ is finite dimensional by a classical theorem of Harder. When $n=0$ the claim is obvious. Assume now that it is true for $n-1$ and pick a point $x\in X$ such that not every element of $W$ is zero at $x$. Then the $\Bbb C$-subspace $V$ of elements of $W$ vanishing at $x$ has dimension $n-1$ hence there is a set $R\subset X-\{x\}$ such that every element of $V$ is uniquely determined by its restriction to $R$. The set $S=R\cup\{x\}$ clearly satisfies the required property.\ $\square$
\enddefinition

\heading 4. Drinfeld modular forms and quotients of the modular jacobian
\endheading

\definition{Definition 4.1} For every proper ideal $\goth m\triangleleft\bold A$ there is a $\goth m$-th Hecke correspondence on the Drinfeld modular curve $X_0(\goth n)$ which in turn induces an endomorphism of the Jacobian $J_0(\goth n)$ of the curve, called the Hecke operator $T_{\goth m}$ (for a detailed description see for example [6] or [7].) The operator $T_{\goth m}$ is denoted by the same symbol we use for the operators introduced in Definition 3.4, but this will not cause confusion as we will see. For the moment it is sufficient to remark that they act on different objects.
Let $\Bbb T'(\goth n)$ denote the algebra with unity generated by the endomorphisms $T_{\goth q}$ of the Jacobian $J_0(\goth n)$, where $\goth q\triangleleft\bold A$ is any prime ideal which does not divide $\goth n$. Let $J^{old}_0(\goth n)$ be the smallest abelian sub-variety of $J_0(\goth n)$ which is left invariant under the action of $\Bbb T'(\goth n)$ and contains the image of $J_0(\goth m)$ with respect to the map $J_0(\goth m)\rightarrow J_0(\goth n)$ induced by the degeneracy map $X_0(\goth n)\rightarrow X_0(\goth m)$ via Picard functoriality for every ideal $\goth m\supsetneqq\goth n$ of $\bold A$. Let $J^{new}_0(\goth n)$ denote the quotient of $J_0(\goth n)$ by $J^{old}_0(\goth n)$. The abelian variety $J^{new}_0(\goth n)$ is naturally equipped with an action of $\Bbb T'(\goth n)$ which makes the quotient map $J_0(\goth n)\rightarrow J^{new}_0(\goth n)$ a $\Bbb T'(\goth n)$-equivariant homomorphism. Let $\Bbb T(\goth n)$ denote the image of $\Bbb T'(\goth n)$ in $\text{End}_F(J^{new}_0(\goth n))$ corresponding to this action and let $T_{\goth q}$ denote the image of the Hecke operator $T_{\goth q}$ under the quotient map $\Bbb T'(\goth n)\rightarrow\Bbb T(\goth n)$ for every prime ideal $\goth q\triangleleft\bold A$ which does not divide $\goth n$ by the usual abuse of notation. The algebra $\Bbb T(\goth n)$ is known to be commutative. Let $\goth E(\goth n)$ denote the ideal of $\Bbb T(\goth n)$ generated by the elements $T_{\goth q}-q^{\deg(\goth q)}-1$, where $\goth q\!\!\not\!|\goth n$ is any prime. The algebra $\Bbb T(\goth n)$ will be called the Hecke algebra and $\goth E(\goth n)$ is its Eisenstein ideal, although these differ slightly from the usual definition, since they do not involve the Atkin-Lehmer operator. The latter will play no role in what follows. By Lemma 4.2 below the $\Bbb T(\goth n)/\goth E(\goth n)$ quotient is finite: we say that a prime
number $l$ is an Eisenstein prime for $\goth n$ if $l$ divides $|\Bbb T(\goth n)/\goth E(\goth n)|$. Note that $J^{old}_0(\goth n)$ is a point when $\goth n$ is a proper prime ideal as $X_0(1)$ is a rational curve. In particular our definition agrees with the one introduced in Definition 7.10 of [19] on page 160 in this particular case.
\enddefinition
\proclaim{Lemma 4.2} The quotient $\Bbb T(\goth n)/\goth E(\goth n)$ is finite and the natural injection $\Bbb Z\rightarrow\Bbb T(\goth n)$ induces a surjective homomorphism $\Bbb Z\rightarrow\Bbb T(\goth n)/\goth E(\goth n)$.
\endproclaim
\definition{Proof} It is clear from the definition that every generator $T_{\goth q}$ of $\Bbb T(\goth n)$ is congruent to an element of $\Bbb Z$ modulo the Eisenstein ideal, so the natural inclusion of $\Bbb Z$ in $\Bbb T(\goth n)$ induces a surjection $\Bbb Z@>>>\Bbb T(\goth n)/\goth E(\goth n)$. If the quotient $\Bbb T(\goth n)/\goth E(\goth n)$ is not finite then this map is also injective, hence the Eisenstein ideal
generates a non-trivial ideal in $\Bbb T(\goth n)\otimes_{\Bbb Z}\Bbb
Q$. Fix a prime $l$ different from the characteristic $p$. Let $\overline F$ denote the separable closure of $F$ and let $J_0(\goth n)_{\overline F}$ denote the base change of $J_0(\goth n)$ to the spectrum of $\overline F$. The action of $\Bbb T'(\goth n)$ on $J_0(\goth n)$ induces an action on the \'etale cohomology group $H^1(J_0(\goth n)_{\overline F},\Bbb Q_l)$ which is faithful. If the quotient $\Bbb T(\goth n)/\goth E(\goth n)$ is not finite then for every prime ideal $\goth q\triangleleft\bold A$ which does not divide $\goth n$ the element $T_{\goth p}-q^{\deg(\goth p)}-1$ generates a non-trivial ideal in the $\Bbb Q_l$-algebra generated by $T_{\goth q}$ considered as an endomorphism of the vector space $H^1(J_0(\goth n)_{\overline F},\Bbb Q_l)$ by the above. Hence $T_{\goth p}-q^{\deg(\goth p)}-1$ divides the minimal polynomial of the endomorphism $T_{\goth q}$. In particular $q^{\deg(\goth p)}+1$ is an eigenvalue of this operator. Let $\Cal J_0(\goth n)$ denote the N\'eron model of $J_0(\goth n)$ over $X$. It is known that $\Cal J_0(\goth n)$ has good reduction over all primes $\goth q$ which does not divide $\goth n$. Let $\text{Frob}_{\goth q}$ denote the Frobenius endomorphism of the  fiber of $\Cal J_0(\goth n)$ over the closed point $\goth q$. By the N\'eron property the Hecke operator $T_{\goth q}$ induces an endomorphism of $\Cal J_0(\goth n)$. It is known that the restriction of the latter to the fiber over $\goth q$ satisfies the Eichler-Shimura relation:
$$\text{Frob}_{\goth q}^2-T_{\goth q}\cdot\text{Frob}_{\goth q}+
q^{\deg(\goth q)}=0.$$
These endomorphisms induce $\Bbb Q_l$-linear operators on the \'etale cohomology group $H^1(\Cal J_0(\goth n)_{\goth q},\Bbb Q_l)$ where $\Cal J_0(\goth n)_{\goth q}$ denotes the base change of the fiber of $\Cal J_0(\goth n)$ over the closed point $\goth q$ to the algebraic closure of the field of definition. By the proper base change theorem there is a $T_{\goth q}$-equivariant isomorphism between the vector spaces $H^1(\Cal J_0(\goth n)_{\goth q},\Bbb Q_l)$ and $H^1(J_0(\goth n)_{\overline F},\Bbb Q_l)$ hence the action of $T_{\goth q}$ on $H^1(\Cal J_0(\goth n)_{\goth q},\Bbb Q_l)$ has eigenvalue $q^{\deg(\goth p)}+1$. Therefore either $1$ or $q^{\deg(\goth q)}$ is an eigenvalue of the operator corresponding to $\text{Frob}_{\goth q}$ by the Eichler-Shimura relation above. But the latter is impossible by Weil's purity theorem if $\deg(\goth q)$ is sufficiently large.\ $\square$
\enddefinition
\definition{Definition 4.3} Recall that for every pair of integers $k$, $m$ a holomorphic function $f:\Omega\rightarrow\Bbb C_{\infty}$ is called a Drinfeld modular form of weight $k$, of type $m$ and of level $\goth n$ if we have
$$f(\gamma z)=(\det\gamma)^{-m}(cz+d)^kf(z)\quad(\forall
z\in\Omega,\forall\gamma=\left(\matrix a&b\\ c&d\endmatrix\right)
\in\Gamma_0(\goth n)),$$
and $f$ is holomorphic at the cusps. (For an explanation of the latter
condition see 2.8.3 of [12], page 46). Also recall that a Drinfeld modular
is called double-cuspidal if its order of vanishing at every cusp is at least two (see 2.8.5 of [12], page 46). Let $M^2_{2,1}(\Gamma_0(\goth n))$ denote the space of double-cuspidal, weight two Drinfeld modular forms of type $1$ and
level $\goth n$. Let $X_{\Bbb C_{\infty}}$ denote the base change of any algebraic variety $X$ over $F$ to the spectrum of $\Bbb C_{\infty}$. For every differential form $\omega\in\Gamma(X_0(\goth n)_{\Bbb C_{\infty}},\Omega^1)$ the pull-back of the restriction $\omega|_{Y_0(\goth n)_{\Bbb C_{\infty}}}$ via the uniformization map of Theorem 3.6 can be written as $f_{\omega}(z)dz$ where $f_{\omega}$ is an element of $M^2_{2,1}(\Gamma_0(\goth n))$. According to Proposition 2.10.2 of [12] on page 47 the map $\Gamma(X_0(\goth n)_{\Bbb C_{\infty}},\Omega^1)\rightarrow M^2_{2,1}(\Gamma_0(\goth n))$ given by the rule $\omega\mapsto f_{\omega}$ is an isomorphism of vector spaces over $\Bbb C_{\infty}$. The vector space $\Gamma(X_0(\goth n)_{\Bbb C_{\infty}},\Omega^1)$ is canonically isomorphic to the $\Bbb C_{\infty}$-dual of the tangent space of $J_0(\goth n)_{\Bbb C_{\infty}}$ at the identity of this group scheme hence it is equipped with an action of the algebra $\Bbb T'(\goth n)$. The corresponding action on $M^2_{2,1}(\Gamma_0(\goth n))$ has a description similar to the one in Definition 4.4.
\enddefinition
\proclaim{Theorem 4.4} There is an $\Bbb C_{\infty}$-vectorspace isomorphism:
$$\text{\rm res}:M^2_{2,1}(\Gamma_0(\goth n))@>>>H_{!!}(\Cal T,\Bbb C_{\infty})^{\Gamma_0(\goth n)}$$
such that $\text{\rm res}(T_{\goth q}(f))=
T_{\goth q}(\text{\rm res}(f))$ for every $f\in M^2_{2,1}(\Gamma_0(\goth n))$ and proper prime ideal $q\triangleleft\bold A$ not dividing $\goth n$.
\endproclaim
\definition{Proof} This is Theorem 6.5.3 of [7] on page 75.\ $\square$
\enddefinition
The aim of this section is to prove the following
\proclaim{Theorem 4.5} The number $p$ is not an Eisenstein prime for $\goth n$ when $\goth n$ is a square-free ideal.
\endproclaim
According to claim $(vi)$ of Proposition 7.11, pages 160-161 and (the proof of) Corollary 11.8, pages 194-195 of [19], the number $|\Bbb T(\goth p)/\goth E(\goth p)|$ is the product of $N(\goth p)$ and a power of $p$ when $\goth p$ is a prime ideal. Hence the result above implies Theorem 1.2. Let $M^{2,new}_{2,1}(\Gamma_0(\goth n))$ denote the subspace of $M^2_{2,1}(\Gamma_0(\goth n))$ corresponding to the
$\Bbb C_{\infty}$-dual of the tangent space of $J^{new}_0(\goth n)_{\Bbb C_{\infty}}$, which is a subspace of the $\Bbb C_{\infty}$-dual of the tangent space of $J_0(\goth n)_{\Bbb C_{\infty}}$, under the identification in Definition 4.3 and let $H^{new}_{!!}(\Cal T,\Bbb C_{\infty})^{\Gamma_0(\goth n)}$ denote its image in $H_{!!}(\Cal T,\Bbb C_{\infty})^{\Gamma_0(\goth n)}$ with respect to the map res of Theorem 4.4. The vector space $H^{new}_{!!}(\Cal T,\Bbb C_{\infty})^{\Gamma_0(\goth n)}$ has a purely combinatorial description as the image of certain new forms under the canonical map $H_!(\Cal T,\Bbb Z)^{\Gamma_0(\goth n)}\otimes\Bbb C_{\infty}\rightarrow H_!(\Cal T,\Bbb C)^{\Gamma_0(\goth n)}$ whose details are left to the dedicated reader. According to the discussion above Theorem 4.5 also has the following immediate
\proclaim{Corollary 4.6} There is no non-zero element of the vector space $M^{2,new}_{2,1}(\Gamma_0(\goth n))$ or $H^{new}_{!!}(\Cal T,\Bbb C_{\infty})^{\Gamma_0(\goth n)}$ fixed by the Hecke operator $T_{\goth q}$ for every proper prime ideal $q\triangleleft\bold A$ not dividing $\goth n$.\ $\square$
\endproclaim
The claim about $H^{new}_{!!}(\Cal T,\Bbb C_{\infty})^{\Gamma_0(\goth n)}$ in the corollary above is of elementary nature about an essentially combinatorial object but its proof is not at all elementary nor combinatorial.
\definition{Definition 4.7} Let $R$ be an order of an algebraic number field $K$. For every abelian variety $A$ defined over a field $\bold k$ and for every prime number $l$ different from the characteristic of $\bold k$ let $T_l(A)$ denote the $l$-th Tate module of $A$ and let $V_l(A)$ denote the vector space $T_l(A)\otimes_{\Bbb Z_l}\Bbb Q_l$. We say that $A$ is equipped with $R$-multiplication $\phi$ over the field $\bold k$ if $\phi:R\rightarrow\text{End}_{\bold k}(A)$ is a non-zero ring homomorphism. In this case the induced $R$-action on $T_l(A)$ makes $V_l(A)$ a free module over $K\otimes_{\Bbb Q}\Bbb Q_l=R\otimes_{\Bbb Z}\Bbb Z_l$ of finite rank. If $\bold k$ is a finite field then the action of the (arithmetic) Frobenius on $V_l(A)$ is $K\otimes_{\Bbb Q}\Bbb Q_l$-linear, the characteristic polynomial of this action lies in $K[t]\subset K\otimes_{\Bbb Q}\Bbb Q_l[t]$ and it is independent of the choice of $l$. Assume now that $A$ is defined over $F$ and let $\Cal A$ and $\Cal A_x$ be the N\'eron model of $A$ over $\Bbb P^1_{\Bbb F_q}$ and the fiber of $\Cal A$ over any closed point $x$ of $\Bbb P^1_{\Bbb F_q}$, respectively. Let $a_x(A)\in K$ denote coefficient of $t$ in the characteristic polynomial $(t)\in K[t]$ of the Frobenius on $V_l(\Cal A_x)$ for every place $x$ of $F$ where $A$ has good reduction. Proper prime ideals of $\bold A$ and the places of $F$ different from $\infty$ are in a natural bijective correspondence. These two sets will be identified in all that follows.
\enddefinition
Assume that $\pi:\Bbb T(\goth n)\rightarrow R$ is a surjective ring homomorphism onto the order $R$.
\proclaim{Theorem 4.8} There is an abelian variety $A$ over $F$ equipped with $R$-multiplication $\phi$ over $F$ such that
\roster
\item"$(i)$" the abelian variety $A$ is the quotient of $J_0^{new}(\goth n)$,
\item"$(ii)$" the dimension of $A$ is equal to the rank of $R$ as a free $\Bbb Z$-module,
\item"$(iii)$" we have $a_{\goth q}(A)=\pi(T_{\goth q})$ for every prime ideal $\goth q\triangleleft\bold A$ which does not divide $\goth n$.
\endroster
\endproclaim
Note that $A$ has good reduction at every prime ideal $\goth q\triangleleft\bold A$ which does not divide $\goth n$ by condition $(i)$ hence condition $(iii)$ is meaningful.
\definition{Proof} This result is just a reformulation of the Langlands correspondence proved in the classical paper [4]. For a detailed discussion see [12].\ $\square$
\enddefinition
\definition{Proof of Theorem 4.5} Assume that $p$ is an Eisenstein prime and let $\goth P\triangleleft\Bbb T(\goth n)$ be a proper minimal prime ideal which is contained in the ideal $(\goth E(\goth n),p)$. Let $R$ denote $\Bbb T(\goth n)/\goth P$ and let $\pi:\Bbb T(\goth n)\rightarrow R$ be the quotient map. Then $R$ is an order in a number field hence there is an abelian variety $A$ over $F$ with the properties described in Theorem 4.8. The image of the ideal $(\goth E(\goth n),p)$ with respect to $\pi$ is a maximal ideal $\goth m$ in $R$ such that $R/\goth m=\Bbb F_p$. For every prime ideal $\goth q\triangleleft\bold A$ which does not divide $\goth n$ we have
$$a_{\goth q}(A)=\pi(T_{\goth q})\equiv q^{\deg(\goth q)}+1\equiv1\mod\goth m$$
by condition $(iii)$. Therefore $A$ has ordinary good reduction at every prime ideal $\goth q\triangleleft\bold A$ which does not divide $\goth n$ by Lemma 4.9 below. On the other hand $J_0^{new}(\goth n)$ has multiplicative reduction at $\infty$ and every prime ideal $\goth q$ dividing $\goth n$ hence so does the abelian variety $A$. But it should have good reduction everywhere by Theorem 2.7 which is a contradiction.\ $\square$
\enddefinition
\proclaim{Lemma 4.9} Let $R$ be an order in a number field $K$, let $\goth m\triangleleft R$ be an ideal such that $R/\goth m=\Bbb F_p$ and let $A$ be an abelian variety equipped with $R$-multiplication $\phi$ over a finite field $\bold k$ of characteristic $p$ such that the dimension of $A$ is equal to the $\Bbb Z$-rank of $R$.
If the coefficient of $t$ in the characteristic polynomial $f(t)\in K[t]$ of the Frobenius on $V_l(A)$ lies in $R-\goth m$ then $A$ is ordinary.
\endproclaim
\definition{Proof} The abelian variety $A$ is ordinary if and only if the multiplicity of the slope $0$ of the Dieudonn\'e-crystal $\Bbb M$ attached to the $p$-divisible group of $A$ is equal to the dimension of $A$. Let $L$ denote the field of fractions of the Witt vectors $\Bbb W(\bold k)$ of $\bold k$ of infinite length. The $R$-multiplication $\phi$ induces an $R$-action on $\Bbb M$ which makes $\Bbb M\otimes_{\Bbb W(\bold k)}L$ a free module over $K\otimes_{\Bbb Q}L=R\otimes_{\Bbb Z}\Bbb W(\bold k)$ of rank two. The action of the $L$-linear Frobenius on $\Bbb M\otimes_{\Bbb W(\bold k)}L$ is $K\otimes_{\Bbb Q}L$-linear and the characteristic polynomial of this action is $f(t)$ by the Weil conjectures. By assumption exactly one of the two roots of $f(t)$ has valuation $0$ as an element of the separable closure of $L$. The $L$-dimension of the free $K\otimes_{\Bbb Q}L$-eigenspace of this eigenvalue is equal to the dimension of $A$ since the latter is equal to the dimension of $K$ as a $\Bbb Q$-vectorspace by assumption. This subspace is left invariant by the absolute Frobenius and its only slope is zero. The claim is now clear.\ $\square$
\enddefinition

\heading 5. Modular symbols
\endheading

\definition{Definition 5.1} A path $\gamma$ on an oriented
graph $G$ is a sequence of edges $$\{\ldots,e_1,e_2,\ldots,e_n,\ldots\}\in\Cal E(G)$$
indexed by the set $I$ where $I=\Bbb Z$, $I=\Bbb N$ or $I=\{0,1,\ldots,m\}$ for some $m\in\Bbb N$ such that $t(e_i)=o(e_{i+1})$ for every $i$, $i+1\in I$. We say that $\gamma$ is an infinite path, a half-line or a finite path whether we are in the first, in the second or in the third case, respectively. For each edge $e\in\Cal E(G)$ let $i_e:\Cal E(G)\rightarrow\Bbb Z$ denote the unique function such that
$$i_e(f)=
\cases+1,&\text{if $f=e$,}\\
-1,&\text{if $f=\overline e$,}\\
\quad\!\!0,&\text{otherwise.}\endcases$$
Let $\gamma$ be a path $\{\ldots,e_1,\ldots,e_n,\ldots\}$ on $G$ such that every edge in $\Cal E(G)$ is only listed finitely many times in the sequence above. Then the function $i_{\gamma}=\sum_{j\in\Bbb Z}i_{e_j}$ is well-defined as the sum above has only finitely many terms non-zero on $e$ for every edge $e\in\Cal E(G)$. Let us consider now the special case $G=\Gamma\backslash\Cal T$ where $\Gamma=\Gamma_0(\goth n)$ is a short-hand notation introduced for convenience. Let $z(\Gamma)$ denote the cardinality of the center of $\Gamma$ and $\Gamma_e$ is the stabilizer of the edge $e\in\Cal E(\Cal T)$ in $\Gamma$. (It is well-known that the latter is finite.) For every path $\gamma$ on the graph $\Gamma\backslash\Cal T$ such that $i_{\gamma}$ is defined in the sense above we define the function $\gamma^*:\Cal E(\Cal T)\rightarrow\Bbb Z$ given by the rule $\gamma^*(e)=|\Gamma_e|i_{\gamma}(\widetilde e)/z(\Gamma)$, where $\widetilde e$ is the image of the edge $e$ in $\Cal E(\Gamma\backslash\Cal T)$ and the absolute sign $|\cdot|$ denotes the cardinality of every finite set. (Since the center of $\Gamma$ leaves the Bruhat-Tits tree invariant, it lies in the stabilizer $\Gamma_e$, therefore the expression above is indeed an integer.)
\enddefinition
\definition{Definition 5.2} Next we are going to define the fundamental arch connecting two different points $a$, $b\in\Bbb P^1(F_{\infty})$ on the Bruhat-Tits tree. We say that a path $\{\ldots,e_1,\ldots,e_n,\ldots\}$ indexed by the set $I$ on an oriented graph $G$ is without backtracking if $\overline{e_i}\neq e_{i+1}$ for every $i$, $i+1\in I$. Let $S(a,b)$ denote the set of those edges of $\Cal T$ which can be represented by a matrix $\left(\smallmatrix\alpha&\beta\\\gamma&\delta\endsmallmatrix\right)$ such that the homogeneous coordinates $(\alpha:\gamma)=a$ and $(\beta:\delta)=b$. The elements of the set $S(a,b)$ can be indexed uniquely by the set of integers such that it becomes an infinite path without backtracking: this is the fundamental arch $\overline{ab}$ connecting $a$ and $b$. Let $\overline{ab}$ denote the image of the fundamental arch under the canonical map $\Cal T\rightarrow\Gamma_0(\goth n)\backslash\Cal T$ as well by slight abuse of notation. Let $[a,b]:\Cal E(\Cal T)\rightarrow\Bbb Z$ denote the function $(\overline{ab})^*$ introduced in Definition 5.1 if the latter is well-defined.
\enddefinition
Let $\goth p\triangleleft\bold A$ be now a proper non-zero prime ideal.
\proclaim{Proposition 5.3} The following holds:
\roster
\item"$(i)$" for every different $a$, $b\in\Bbb P^1(F)$ the function $[a,b]$ is well-defined and it is a $\Bbb Z$-valued left ${\Gamma_0(\goth p)}$-invariant harmonic cochain,
\item"$(ii)$" we have $[a,b]\in H_!(\Cal T,\Bbb Z)^{\Gamma_0(\goth p)}$ for every different $a$, $b\in\Bbb P^1(F)$ which are equivalent under the M\"obius action of $\Gamma_0(\goth p)$,
\item"$(iii)$" for every proper non-zero prime ideal $\goth p\neq\goth q\triangleleft \bold A$ we have
$$(1+q^{\deg(\goth q)}-T_{\goth q})[0,\infty]\in
(q-1)H_!(\Cal T,\Bbb Z)^{\Gamma_0(\goth p)}.$$
\endroster
\endproclaim
\definition{Proof} We say that two half-lines on an oriented graph are equivalent if they only differ in a finite graph. We also say that a half-line on $\Cal T$ is $F$-rational if it is without backtracking and it is contained in a fundamental arch $\overline{ab}$ for some $a$, $b\in\Bbb P^1(F)$. By Serre's structure theorem the graph $\Gamma_0(\goth p)\backslash\Cal T$ is the union of a finite graph and finitely many half-lines. Moreover every $F$-rational half-line on $\Cal T$ is equivalent to a half-line $\gamma$ such that the restriction of the projection $\Cal E(\Cal T)\rightarrow\Cal E(\Gamma_0(\goth p)\backslash\Cal T)$ onto $\gamma$ is a bijection onto one of those half-lines above. In particular the functions introduced in Definition 5.1 are all well-defined for the image of every $F$-rational half-line under the projection above. For every $a$, $b\in\Bbb P^1(F)$ there are two half-lines $\alpha$ and $\beta$ on $\Cal T$ such that the disjoint union of $\alpha$ and the set we get from $\beta$ by reversing each of its edges is equal to $\overline{ab}$ as a set. Let $\alpha$ and $\beta$ also denote the image of $\alpha$ and $\beta$ under the canonical map $\Cal T\rightarrow\Gamma_0(\goth p)\backslash\Cal T$, respectively. Then $i_{\overline{ab}}=i_{\alpha}-i_{\beta}$ as functions on $\Cal E(\Gamma_0(\goth p)\backslash\Cal T)$. Hence the function $[a,b]$ is well-defined. Since $\gamma^*$ is a harmonic cochain for every infinite path $\gamma$ on $\Gamma_0(\goth p)\backslash\Cal T$ whenever the former is defined, we get that claim $(i)$ is true. If $a$ and $b$ are equivalent under the M\"obius action of $\Gamma_0(\goth p)$ then the half-lines $\alpha$ and $\beta$ on $\Gamma_0(\goth p)\backslash\Cal T$ are equivalent hence the difference $i_{\alpha}-i_{\beta}$ is supported on a finite set. Therefore claim $(ii)$ is also true. We continue with the proof of part $(iii)$. Let  $r$ be the unique monic polynomial generating $\goth q$. Let $R(\goth q)\subset\Bbb F_q[T]$ denote the set of non-zero polynomials whose degree is less than $\deg(\goth q)$ and let $R(\goth q,\goth p)$ denote the set
$$R(\goth q,\goth p)=\{\left(\matrix1&0\\0&r\endmatrix\right)\}
\cup\{\left(\matrix r&a\\0&1\endmatrix\right)|a\in R(\goth q)\}.$$
Then the set $R(\goth q,\goth p)$ is a set of representatives of the double coset $H(\goth q,\goth p)$ introduced in Definition 3.4 as the notation indicates. For every $e\in\Cal E(\Cal T)$ we have:
$$\split T_{\goth q}([0,\infty])(e)&={1\over z(\Gamma)}\!\!
\sum_{M\in R(\goth q,\goth p)}\!\!\!\!\left(
|\{\gamma\in\Gamma|e\in M\gamma[0,\infty]\}|-
|\{\gamma\in\Gamma|\overline e\in M\gamma[0,\infty]\}|\right)\\
&={1\over z(\Gamma)}\!\!\sum_{M\in R(\goth q,\goth p)}
\!\!\!\!\left(|\{\gamma\in\Gamma|e\in \gamma M[0,\infty]\}|-
|\{\gamma\in\Gamma|\overline e\in \gamma M[0,\infty]\}|\right)\\
&=2[0,\infty]+\sum_{0\neq a\in R(\goth q)}[a/r,\infty]\endsplit$$
using again the notation $\Gamma=\Gamma_0(\goth p)$ of Definition 5.1. The first equation above is just the definition. The second equation follows from the fact that every product $M\gamma$ where $M\in R(\goth q)$ and $\gamma\in\Gamma_0(\goth p)$ can be written uniquely as a product $\gamma'M'$ for some $M'\in R(\goth q)$ and $\gamma'\in\Gamma_0(\goth p)$. The last equation follows from the fact that for every $M=\left(\smallmatrix\alpha&\beta\\\gamma&\delta\endsmallmatrix\right)\in GL_2(F)$ the image of the fundamental arch $\overline{0\infty}$ under the automorphism of $\Cal T$ induced by multiplication on the left by the matrix $M$ is the fundamental arch $\overline{(\beta:\delta)(\alpha:\gamma)}$.
Note that for every $a$, $b$ and $c\in\Bbb P^1(F)$ which are pair-wise different we have $[a,b]=[a,c]+[c,b]$. Therefore
$$(1+q^{\deg(\goth q)}-T_{\goth q})[0,\infty]=\sum_{0\neq a\in R(\goth q)}[0,a/r].$$
Note that for every $0\neq a\in R(\goth q)$ the polynomials $r$ and $a$ are relatively prime hence there are polynomials $c$, $d\in\bold A$ such that the matrix $\left(\smallmatrix c&a\\ d&r\endsmallmatrix\right)$ is an element of $\Gamma_0(\goth p)$. This matrix maps $0$ to $a/r$ via the M\"obius action so $[0,a/r]\in H_!(\Cal T,\Bbb Z)^{\Gamma_0(\goth p)}$ by claim $(ii)$. Since $[0,a/r]=[0,ca/r]$ for every $a\in R(\goth q)$ and every $c\in\Bbb F^*_q$ the equation above implies claim $(iii)$.\ $\square$
\enddefinition
\definition{Definition 5.4} Recall that a finite path $\{e_0,e_2,\ldots,e_n\}\in\Cal E(G)$ on an oriented graph $G$ is closed if the equality $t(e_n)=o(e_0)$ holds, too. We define $H_1(G,\Bbb Z)$ as the abelian group of $\Bbb Z$-valued functions on $\Cal E(G)$ generated by the functions $i_{\gamma}$ where $\gamma$ is a closed path. We define the map
$$j_{\Gamma_0(\goth n)}:H_1(\Gamma_0(\goth n)\backslash\Cal T,\Bbb Z)\rightarrow H_!(\Cal T,\Bbb Z)^{\Gamma_0(\goth n)},$$
as the unique homomorphism which maps $i_{\gamma}$ to the cochain $\gamma^*$ for every $\gamma$ closed path, using the notations of Definition 5.1. It is easy to see that the homomorphism is
well-defined, that is $\gamma^*$ is indeed a harmonic cochain. By a theorem of Gekeler and Nonnengardt (Theorem 3.3 of [11], page 702) this homomorphism is in fact an isomorphism.
\enddefinition
\definition{Definition 5.5} Let
$\Gamma_0(\goth p)_{\text{ab}}=\Gamma_0(\goth p)/[\Gamma_0(\goth p),\Gamma_0(\goth p)]$ be the abelianization of
$\Gamma_0(\goth p)$, and let $\overline{\Gamma}_0(\goth p)=\Gamma_0(\goth p)_{\text{ab}}/
(\Gamma_0(\goth p)_{\text{ab}})_{\text{tors}}$ be its maximal torsion-free quotient. For each $\gamma\in\Gamma_0(\goth p)$ let $\overline{\gamma}$ denote its image in $\overline{\Gamma}_0(\goth p)$. Fix a vertex $v_0\in\Cal V(\Cal T)$ and for every $\gamma\in\Gamma_0(\goth p)$ let $e_0,e_1,\ldots,e_{n(\gamma)}$ be the unique geodesic path connecting $v_0$ with $\gamma(v_0)$, that is $v_0=o(e_1)$ and $\gamma(v_0)=t(e_{n(\gamma)})$. Recall that a path is geodesic if it is the shortest connecting its endpoints, in this case $v_0$ with $\gamma(v_0)$, i.e$\text{.}$ the number $n(\gamma)$ is the smallest possible. The image of the path $e_0,e_1,\ldots$ is closed in $\Gamma_0(\goth p)
\backslash\Cal T$: let $i(\gamma)$ denote the corresponding element in $H_1(\Gamma_0(\goth p)\backslash\Cal T,\Bbb Z)$. The function $i$ induces a homomorphism $i:\overline{\Gamma}_0(\goth p)\rightarrow H_1(\Gamma_0(\goth p)\backslash\Cal T,\Bbb Z)$ which is independent of the choice of $v_0$ and it is an
isomorphism. We will use this identification without further notice.
\enddefinition
Let $\Phi_{AJ}:\text{\rm Hom}(\overline{\Gamma}_0(\goth p),\Bbb C_{\infty}^*)
\rightarrow J_0(\goth p)(\Bbb C_{\infty})$ be the Abel-Jacobi map of Gekeler-Reversat. (For a description using the same notation see sections 7.1-7.7 of [19] on pages 158-159.) It is a group homomorphism holomorphic in the rigid analytic sense. Let $c:
\overline{\Gamma}_0(\goth p)\rightarrow\text{\rm Hom}
(\overline{\Gamma}_0(\goth p),\Bbb C_{\infty}^*)$ be the period map defined in Proposition 7.5 of [19] on page 159. Note that the domain of the map $c$ is equipped with an action of the operator $T_{\goth q}$ via the isomorphism $j_{\Gamma_0(\goth p)}$. The following result is Theorem 7.9 of [19] on page 139.
\proclaim{Theorem 5.6} For every prime $\goth q\triangleleft\bold A$, different
from $\goth p$, there is a unique endomorphism $T_{\goth q}$ of the rigid
analytic torus $\text{\rm Hom}(\overline{\Gamma}_0(\goth p),
\Bbb C_{\infty}^*)$, which leaves the lattice $\overline{\Gamma}_0(\goth
p)$ invariant, and makes the diagram:
$$\CD0@>>>\overline{\Gamma}_0(\goth p)@>c>>\text{\rm Hom}
(\overline{\Gamma}_0(\goth p),\Bbb C_{\infty}^*)@>\Phi_{AJ}>>J_0(\goth
p)(\Bbb C_{\infty})@>>>0\\ @.@VT_{\goth q}VV@VT_{\goth q}VV@VT_{\goth q}VV\\
0@>>>\overline{\Gamma}_0(\goth p)@>c>>\text{\rm Hom}(\overline{\Gamma}_0
(\goth p),\Bbb C_{\infty}^*)@>\Phi_{AJ}>>J_0(\goth p)(\Bbb C_{\infty})@>>>0\endCD$$
commutative. Moreover the map $j:\overline{\Gamma}_0(\goth p)@>>>\Cal
H_0(\goth p,\Bbb Z)$ is equivariant with respect to this action on
$\overline{\Gamma}_0(\goth p)$ and the action of the Hecke operator
$T_{\goth q}$ on $\Cal H_0(\goth p,\Bbb Z)$.\ $\square$
\endproclaim
\definition{Remarks 5.7} For the moment let $\Bbb T(\goth p)'$ denote the ring with unity generated by the operators $T_{\goth q}$ acting on the torus $\text{\rm Hom}(\overline{\Gamma}_0(\goth p),\Bbb C_{\infty}^*)$. The rigid analytic endomorphisms of algebraic tori are algebraic, so they act faithfully on any Zariski-dense invariant subset. In particular the action of $\Bbb T(\goth p)'$ on the module $\overline{\Gamma}_0(\goth p)$ is faithful.  Since $\Phi_{AJ}$ injects $\text{\rm Hom}(\overline{\Gamma}_0(\goth p),\Cal O_{\infty}^*)$ into $J_0(\goth p)(F_{\infty})$, the action of $\Bbb T(\goth p)'$ on the Jacobian $J_0(\goth p)$ is also faithful. Therefore the algebra  $\Bbb T(\goth p)'$ is canonically isomorphic to the Hecke algebra $\Bbb T(\goth p)$ introduced in Definition 4.1. In particular we have a well-defined homomorphism
$$e:\goth E(\goth p)\rightarrow H_!(\Cal T,\Bbb Z)^{\Gamma_0(\goth p)}$$
of $\Bbb T(\goth p)$-modules given by the rule $\alpha\mapsto\alpha([0,\infty])/(q-1)$ according to part $(iii)$ of Proposition 5.3. This map is the analogue of the winding
homomorphism introduced by Mazur.
\enddefinition
\definition{Definition 5.8} Let $B$ denote the group scheme of invertible upper triangular two by two matrices. Let $h:\Gamma_0(\goth p)\rightarrow\Bbb Z/N(\goth p)\Bbb Z$ denote the composition of the reduction map $r:\Gamma_0(\goth p)\rightarrow
 B(\bold A/\goth p)\subset GL_2(\bold A/\goth p)$ mod $\goth p$, the upper left conner element $a:B(\bold A/\goth p)\rightarrow(\bold A/\goth p)^*$ and the surjection $s:(\bold A/\goth p)^*\rightarrow\Bbb Z/N(\goth p)\Bbb Z$ which unique up to isomorphism. This homomorphism factors through $\overline{\Gamma}_0(\goth p)$, so it induces a homomorphism:
$$\phi=h\circ j_{\Gamma_0(\goth n)}^{-1}:H_!(\Cal T,\Bbb Z)^{\Gamma_0(\goth p)}\rightarrow\Bbb Z/N(\goth p)\Bbb Z$$
under the identification of Definition 5.5.
\enddefinition
\proclaim{Lemma 5.9} The kernel of $\phi$ contains $\goth E(\goth p)
H_!(\Cal T,\Bbb Z)^{\Gamma_0(\goth p)}$ and induces an isomorphism
$$\overline{\phi}:H_!(\Cal T,\Bbb Z)^{\Gamma_0(\goth p)}/
\goth E(\goth p)H_!(\Cal T,\Bbb Z)^{\Gamma_0(\goth p)}\rightarrow
\Bbb Z/N(\goth p)\Bbb Z.$$
\endproclaim
\definition{Proof} Let $\mu_{\infty}\subset\Bbb C^*_{\infty}$ denote the subgroup of the roots of unity. Under the map $\Psi_{AJ}$ of Theorem 5.6 the image of Hom$(\overline{\Gamma}_0(\goth p),\mu_{\infty})$ maps isomorphically onto the group of those torsion points of $J_0(\goth p)(\Bbb C_{\infty})$ which are defined over the maximal unramified extension $K_{\infty}$ of $F_{\infty}$ and map into the connected component of the identity of the special fiber of the N\'eron model of $J_0(\goth p)$ over the valuation ring of $K_{\infty}$ under the specialization map. By Theorem 1.2 the prime $p$ does not divide the order of $\Bbb T(\goth p)/\goth E(\goth p)$ hence the quotient $H_!(\Cal T,\Bbb Z)^{\Gamma_0(\goth p)}/\goth E(\goth p)H_!(\Cal T,\Bbb Z)^{\Gamma_0(\goth p)}$ is dual to the part of Hom$(\overline{\Gamma}_0(\goth p),\mu_{\infty})$ annihilated by the Eisenstein ideal. There is a unique subgroup of $\mu_{\infty}$ isomorphic to $\Bbb Z/N(\goth p)\Bbb Z$ hence the homomorphism $\phi$ can be con\-sidered as an element of
Hom$(\overline{\Gamma}_0(\goth p),\mu_{\infty})$. During the proof of claim $(i)$ of Proposition 8.18 of [19], pages 171-172, we also showed incidentally that the subgroup generated by $\phi$ is the Shimura group $\Cal S(\goth p)$ introduced in the paper quoted above. On the other hand some of the main results of this paper, in particular Theorem 10.5 on page 187 and Corollary 11.8 on pages 194-195 imply that the Shimura group is the largest subgroup of Hom$(\overline{\Gamma}_0(\goth p),\mu_{\infty})$ annihilated by the Eisenstein ideal.\ $\square$
\enddefinition
For every element $b\in\bold A$ let $\overline b\in\bold A/\goth p$ denote the reduction of $b$ modulo the ideal $\goth p$. The following proposition is the analogue of Mazur's congruence formula for the modular symbol:
\proclaim{Proposition 5.10} Let $a$, $b$ be relatively prime elements in $\bold A$ with $b$ relatively prime to $\goth q$. Then
$$\phi([0,a/b])=s(\overline b)\in\Bbb Z/N(\goth p)\Bbb Z.$$
\endproclaim
\definition{Proof}  By assumption there are polynomials $c$, $d\in\bold A$ such that the matrix $M=\left(\smallmatrix c&a\\ d&b\endsmallmatrix\right)$ is an element of $\Gamma_0(\goth p)$. This matrix maps $0$ to $a/b$ so $[0,a/b]\in H_!(\Cal T,\Bbb Z)^{\Gamma_0(\goth p)}$ by claim $(ii)$ of Proposition 5.3. In particular the left hand side of the equation above is well-defined.
If $\alpha$ is a half-line on $\Cal T$ such that the set we get from $\alpha$ by reversing each of its edges is contained in $\overline{0(a/b)}$ then the half-line $M\alpha$ is equivalent to a half-line contained in $\overline{0(a/b)}$ because the images of the half-lines $\alpha$ and $M\alpha$ with respect to the canonical map $\Cal T\rightarrow\Gamma_0(\goth n)\backslash\Cal T$
are equivalent. Therefore there is an edge $e\in\overline{0(a/b)}$ such that $\overline{Me}\in\overline{0(a/b)}$, too. Now let $\alpha$ denote the unique half line whose first edge (indexed by zero) is $\overline e$ and $M\alpha$ is contained in $\overline{0(a/b)}$. Let $\gamma$ denote the unique geodesic connecting the vertex $o(e)$ to the vertex $o(Me)=Mo(e)$ and let $\gamma'$ denote its image with respect to the canonical map $\Cal T\rightarrow\Gamma_0(\goth n)\backslash\Cal T$. Then $\overline{0(a/b)}$ as a set is the union of $M\alpha$, $\gamma$ and set we get from $\alpha$ by reversing each of its edges. Therefore $(\overline{0(a/b)})^*=(\gamma')^*$ using the notation of Definition 5.1, so $[0,a/b]=j_{\Gamma(\goth p)}(i(M))$ using the notation of Definitions 5.4 and 5.5. But $\phi(i(M))=s(\overline b)$ by definition.\ $\square$
\enddefinition
For every prime $\goth q\triangleleft\bold A$, different from $\goth p$, let $\overline{\goth q}\in\bold A/\goth p$ denote the reduction $\overline{r(T)}$ of the unique monic polynomial $r(T)\in\bold A$ generating the ideal $\goth q$ modulo the ideal $\goth p$. The next claim is the analogue of Mazur's congruence formula for the
winding homomorphism:
\proclaim{Corollary 5.11} For every prime $\goth q\triangleleft\bold A$, different from $\goth p$, we have
$$\phi\circ e(1+q^{\deg(\goth q)}-
T_{\goth q})={q^{\deg(\goth q)}-1\over q-1}s(\overline{\goth
q})\in\Bbb Z/N(\goth p)\Bbb Z.$$
\endproclaim
\definition{Proof} Let $S(\goth q)\subset\Bbb F_q[T]$ denote the set of non-zero monic polynomials of degree less than deg$(\goth q)$. According to the formula which we derived in the proof of Proposition 5.3 we have:
$$\phi\circ e(1+q^{\deg(\goth q)}-
T_{\goth q})=\phi(\sum_{a\in S(\goth q)}[0,a/r])={q^{\deg(\goth q)}-1\over q-1}s(\overline{\goth
q}),$$
where we used Proposition 5.10 in the second equation.\ $\square$
\enddefinition
\definition{Definition 5.12} For the rest of this chapter we fix an Eisenstein prime $l$ for $\goth p$. We define the $\Bbb Z_l$ algebra $\Bbb T_l(\goth p)$ as the tensor product $\Bbb T(\goth p)\otimes\Bbb Z_l$ and let $\goth P\triangleleft\Bbb T_l(\goth p)$ be the unique prime ideal lying above the ideal $\goth I$ generated by $\goth E(\goth p)$ in $\Bbb T_l(\goth p)$. Let $\Bbb T_{\goth P}$ denote the completion of $\Bbb T_l(\goth p)$ at the prime ideal $\goth P$: this algebra is canonically isomorphic to the completion of $\Bbb T(\goth p)$ at the ideal $\goth P(\goth p,l)$ where $\goth P(\goth p,l)
\triangleleft\Bbb T(\goth p)$ is the unique prime ideal lying above the ideal generated by $\goth E(\goth p)$ and $l$. As $\Bbb Z$ surjects
onto $\Bbb T(\goth p)/\goth E(\goth p)$ via its natural inclusion into $\Bbb T(\goth p)$, clearly $\goth P(\goth p,l)=(\goth E(\goth p),l)$. Hence the latter is a maximal ideal with residue field $\Bbb F_l$.
Let $\eta_{\goth q}$ denote the element $T_{\goth q}-q^{\deg(\goth q)}-1\in\Bbb T(\goth p)$, where $\goth q\triangleleft\bold A$ is any prime ideal different from $\goth p$. Now fix such an ideal $\goth q\triangleleft\bold A$ and let $r(T)\in\bold A$ denote again the unique monic polynomial which generates $\goth q$. We say that $\goth q$ is a good prime if the element
$${q^{\deg(\goth q)}-1\over q-1}s(\overline{\goth
q})\in\Bbb Z/N(\goth p)\Bbb Z$$
is not divisible by $l$ in the group $\Bbb Z/N(\goth p)\Bbb Z$. It is very easy to see that this definition is equivalent to notion we introduced in Definition 10.1 of [19] on page 186. As we already remarked there, the Chebotarev density theorem implies that there are infinitely many good primes. Let $H_{\goth P}$ denote the $\Bbb T_{\goth P}$-module $H_!(\Cal T,\Bbb Z)^{\Gamma_0(\goth p)}\otimes_{\Bbb T(\goth p)}
\Bbb T_{\goth P}$ and let $e_l:\goth I\rightarrow H_{\goth P}$ denote the tensor product $e\otimes_{\Bbb T(\goth p)}$id.
\enddefinition
\proclaim{Theorem 5.13} The following holds:
\roster
\item"$(i)$" the module $H_{\goth P}$ is a free $\Bbb T_{\goth P}$-module of rank one,
\item"$(ii)$" the winding homomorphism $e_l:\goth I\rightarrow H_{\goth P}$ is an isomorphism,
\item"$(iii)$" for every $\goth q\triangleleft\bold A$ prime ideal different from $\goth p$ the element $\eta_{\goth q}$ is a generator of the ideal $\goth I$ if and only if $\goth q$ is a good prime.
\endroster
\endproclaim
\definition{Proof} The module $H_{\goth P}$ is just the base change of the finitely generated $\Bbb T_l(\goth p)$-module $H_!(\Cal T,\Bbb Z_l)^{\Gamma_0(\goth p)}$ to the ring $\Bbb T_{\goth P}$. The latter is the $\Bbb Z_l$-dual of a $\Bbb T_l(\goth p)$-module $T_l$ introduced in Definition 7.10 of [19] on page 160. The module $T_l$ is proved to be a locally free $\Bbb T_l(\goth p)$-module of rank one in Proposition 7.11 of the paper quoted above. Since $\Bbb T_{\goth P}$ is a Gorenstein ring over $\Bbb Z_l$ by Theorems 10.2 and 11.6 of the same paper, the module $H_{\goth P}$ must be free of rank one as the claim $(i)$ says above. Reducing the winding homomorphism $e_l$ mod $\goth I$ one gets a homomorphism $\epsilon:\goth I/\goth I^2\rightarrow H_{\goth P}/\goth I$ and by Lemma 5.9 and Corollary 5.11 the latter maps the element $\eta_{\goth q}$ to a generator of $H_{\goth P}/\goth I$ if and only if $\goth q$ is a good prime. Since there are good primes the map $e_l$ must be surjective by Nakayama's lemma. But $e_l$ is a map between torsion-free modules whose base change to $\Bbb T_l(\goth p)\otimes_{\Bbb Z_l}\Bbb Q_l$ has dimension one, so it must be an isomorphism. Hence $\goth I$ is free and by the above and the element $\eta_{\goth q}$ is a generator if and only if $\goth q$ is a good prime.\ $\square$
\enddefinition

\heading 6. Universal deformation rings
\endheading

\definition{Notation 6.1} For every field $K$ let $\overline K$ and Gal$(\overline K|K)$ denote the separable closure and the absolute Galois group of $K$, respectively. For every prime number $l\neq p$ let $\chi_l:\text{Gal}(\overline F|F)\rightarrow \Bbb Z_l^*$ denote the cyclotomic character and let $\overline{\chi}_l:\text{Gal}(\overline F|F)\rightarrow\Bbb F_l^*$ denote the composition of $\chi_l$ and the reduction $\Bbb Z_l^*\rightarrow\Bbb F_l^*$ mod $l\Bbb Z_l$. Moreover let $\overline{\phi}_l:\text{Gal}(\overline F|F)\rightarrow\Bbb F_l$ denote the unique surjective additive homomorphism unramified at every place of $F$ such that $\overline{\phi}_l(\text{Fr}_{\infty})=1$ where $\text{Fr}_{\infty}$ is a lift of the Frobenius at the place $\infty$ corresponding to the point at infinity on the projective line $\Bbb P^1_{\Bbb F_q}$ over $\Bbb F_q$. Finally let $\overline{\rho}_l:\text{Gal}(\overline F|F)\rightarrow GL_2(\Bbb F_l)$ denote the Galois representation such that
$$\overline{\rho}_l(g)=\left(\matrix1&0\\0&\overline{\chi}_l(g)\endmatrix\right)
\quad(\forall g\in\text{Gal}(\overline F|F))$$
when $l$ does not divide $q-1$ and
$$\overline{\rho}_l(g)=\left(\matrix1&\overline{\phi}_l(g)\\0&1\endmatrix\right)
\quad(\forall g\in\text{Gal}(\overline F|F))$$
when $l$ does divide $q-1$.
\enddefinition
\definition{Definition 6.2} For any commutative ring $A$ with unity and $A$-module $M$ let $GL_A(M)$ denote the group of $A$-linear automorphisms of $M$. For any local ring $A$ let $\goth m_A$ denote the maximal ideal of $A$. Let $\bold k$ be a field and let $\Cal O$ be a noetherian complete local ring with residue field $\bold k$. let $\Cal C$ denote the category of those local Artinian $\Cal O$-algebras such that the map $\Cal O\rightarrow A/\goth m_A$ induced by the $\Cal O$-structure is surjective for every object $A$ of $\Cal C$. Let $G$ be a topological group. A representation of $G$ over an object $A$ of $\Cal C$ is an ordered pair $(M,\rho)$ where $M$ is a finitely generated free $A$-module and $\rho$ is a continuous $G$-action $\rho:G\rightarrow GL_A(M)$ where $M$ is equipped with the discrete topology. For every object $A$ of the category $\Cal C$ let $\Cal P_G(A)$ denote the class consisting of ordered pairs $(V,L)$ where $V$ is a representation $(M,\rho)$ of $G$ over $A$ and $L\subseteq M$ is a $A$-submodule which is a direct summand of $M$ and it is free as an $A$-module. Two elements $O=((M,\rho),V)$ and $O'=((M',\rho'),L')$ of $\Cal P_G(A)$ are said to be isomorphic if there is an isomorphism $\phi:M\rightarrow M'$ of $A[G]$-modules such that $\phi(L)=L'$. Such a relation will be denoted by $O\cong_{A[G]}O'$. For every morphism $f:A\rightarrow A'$ and every element $O=((M,\rho),L)\in\Cal P_G(A)$ let $O\otimes_AA'$ denote the couple $((M\otimes_AA',\rho_{A'}),L\otimes_AA')$ where $\rho_{A'}$ is the composition of $\rho$ and the map $GL_A(M)\rightarrow GL_{A'}(M\otimes_AA')$ which assigns $\phi\otimes_A\text{id}_{A'}$ to $\phi\in GL_A(M)$. Let $O$ be an element of $\Cal P_G(\bold k)$. By a deformation of $O$ over $A\in\text{Ob}(\Cal C)$ we mean an isomorphism class of elements $P$ of $\Cal P_G(A)$ which satisfy the condition $P\otimes_A\bold k\cong_{\bold k[G]}O$. Let Def$(O,A)$ denote the set of deformations of $O$ over $A$. Every morphism $f:A\rightarrow A'$ induces a map $f_*:\text{Def}(O,A)\rightarrow\text{Def}(O,A')$ that maps the isomorphism class of  $O\in\Cal P_G(A)$ to the isomorphism class of $O\otimes_AA'\in\Cal P_G(A')$ which makes Def$(O,\cdot)$ a covariant functor from $\Cal C$ to the category $Sets$ of sets.
\enddefinition
\definition{Definition 6.3} For any commutative ring $A$ with unity let $A[\epsilon]$ denote the $A$-algebra generated by $\epsilon$ subject to the relation $\epsilon^2=0$. When $A\in\text{Ob}(\Cal C)$ so does $A[\epsilon]$. The unique surjective $A$-algebra homomorphism $a:A[\epsilon]\rightarrow A$ is called the augmentation map. For every pair of local $\Cal O$-algebras $A$ and $B$ let $\text{Hom}_{\Cal O}(A,B)$ denote the $\Cal O$-module of local $\Cal O$-algebra homomorphisms. We say that a covariant functor $D:\Cal C\rightarrow Sets$ is pro-representable if there is a complete local $\Cal O$-algebra $R$ such that $R/\goth m_R^n$ is an object of $\Cal C$ for every natural number $n$ and the functors $D$ and $\text{Hom}_{\Cal O}(R,\cdot)$ are naturally isomorphic. Let $D:\Cal C\rightarrow Sets$ be a covariant functor such that $D(\bold k)$ consists of a single set and let
$$
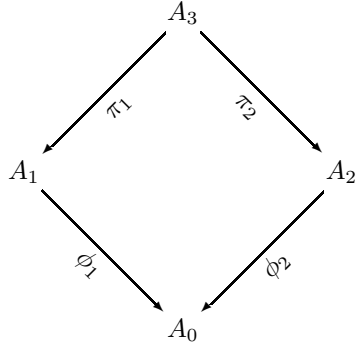
\diagram & &A_3& &\\
&\ldTo_{\pi_1}& &\rdTo_{\pi_2}& \\
A_1& & & &A_2\\
&\rdTo_{\phi_1}& &\ldTo_{\phi_2}& \\
& &A_0& &\\\enddiagram$$
be a cartesian diagram of artinian rings in $\Cal C$ which means that the diagram is commutative and $A_3$ is the fiber product $A_1\times_A A_2$. Because the compositions
$D(\phi_1)\circ D(\pi_1):D(A_3)\rightarrow D(A_0)$ and $D(\phi_2)\circ D(\pi_2):D(A_3)\rightarrow D(A_0)$ are equal we have a map $h:D(A_3)\rightarrow D(A_1)\times_{D(A_0)}D(A_2)$. Assume for a moment that $A_1=A_2=\bold k[\epsilon]$ and $\phi_1=\phi_2$ is the augmentation map. If the morphism
$$h:D(\bold k[\epsilon]\times_{\bold k}\bold k[\epsilon])\rightarrow D(\bold k[\epsilon])\times D(\bold k[\epsilon])$$
is bijective then the set $t_D=D(\bold k[\epsilon])$ is equipped with a natural
$\bold k$-vectorspace structure according to Lemma 2.10 of [21] on page 212 and it is called the tangent space of the functor $D$. We say that a morphism $h:A\rightarrow B$ of $\Cal C$ is small if it is surjective and its kernel is a principal ideal annihilated by $\goth m_A$.
\enddefinition
\proclaim{Theorem (Schlessinger) 6.4} Let $D:\Cal C\rightarrow Sets$ be a covariant functor such that $D(\bold k)$ consists of a single set. Then $D$ is pro-representable by a noetherian local $\Cal O$-algebra if and only if the following conditions hold:
\roster
\item"{\bf H1.}" the map $h$ is surjective if $\phi_1:A_1\rightarrow A_3$ is small,
\item"{\bf H2.}" the map $h$ is bijective if $\phi_2:A_1\rightarrow A_3$ is the augmentation map $a:\bold k[\epsilon]\rightarrow\bold k$,
\item"{\bf H3.}" the dimension $\dim_{\bold k}(t_D)$ is finite,
\item"{\bf H4.}" the map $h$ is bijective if $\phi_1:A_1\rightarrow A_3$ and $\phi_2:A_2\rightarrow A_3$ are equal and small,
\endroster
using the notation of Definition 6.3.
\endproclaim
Note that condition {\bf H3} makes sense because of condition {\bf H2}.
\definition{Proof} This is Theorem 2.11 of [21] on pages 212-215.\ $\square$
\enddefinition
\definition{Definition 6.5} For every place $x$ of $F$ let $D_x\subset\text{Gal}(\overline F|F)$ be the decomposition group at $x$ (unique up to conjugation) and let $I_x\triangleleft D_x$ denote its inertia subgroup. Fix a prime number $l\neq p$ and suppose that $\Cal O=\Bbb Z_l$ and $G=\text{Gal}(\overline F|F)$ using the notation of Definitions 6.2-6.3. Let $O_l\in\Cal P_G(\Bbb F_l)$ denote the couple $((\Bbb F_l^2,\rho_l),\overline L)$ where $\overline L\subset\Bbb F_l^2$ is an affine line which is not fixed by $\text{Gal}(\overline F|F)$ with respect to $\rho_l$. (To fix ideas, let $\overline L$ be spanned by $(1,1)$, when $l\!\not|q-1$, and let $\overline L$ be spanned by $(0,1)$, otherwise.) Let $S$ be a finite subset of the set of proper prime ideals of $\bold A$. For every $A\in\text{Ob}(\Cal C)$ let $\text{Def}_S(A)\subseteq\text{Def}(O_l,A)$ denote the subset of isomorphism classes of those deformations $O=((M,\rho),L)$ of $O_l$ which satisfy the following conditions:
\roster
\item"{\bf D1.}" the representation $\rho$ is unramified at every prime ideal $\goth q\not\in S$,
\item"{\bf D2.}" the inertia subgroup $I_{\goth p}$ at $\goth p$ acts trivially on the submodule $L$,
\item"{\bf D3.}" the determinant of $\rho$ is equal to the composition $\chi_A$ of the cyclotomic character $\chi_l:\text{Gal}(\overline F|F)\rightarrow\Bbb Z_l^*$ and the natural map $\Bbb Z_l^*\rightarrow A^*$,
\item"{\bf D4.}" the $A[D_{\infty}]$-module $M$ has a filtration:
$$\CD0@>>>A(1)@>>>M@>>>A@>>>0\endCD$$
where $A(1)$, $A$ denotes the free $A$-module of rank one equipped with the $D_{\infty}$-action given by $\chi_A|_{D_{\infty}}$ and the trivial character, respectively.
\endroster
For every morphism $f:A\rightarrow A'$ the function $f_*:\text{Def}(O_l,A)\rightarrow\text{Def}(O_l,A')$ maps $\text{Def}_S(A)$ into $\text{Def}_S(A')$ hence $\text{Def}_S(\cdot)$ is a sub-functor of $\text{Def}(O_l,\cdot)$.
\enddefinition
The aim of the chapter is to prove theorem below with the aid of Schlessin\-ger's criteria.
\proclaim{Theorem 6.6} The covariant functor $\text{\rm Def}_S(\cdot)$ is pro-representable by a noetherian complete local $\Cal O$-algebra.
\endproclaim
The noetherian complete local $\Bbb Z_l$-algebra $R$ which pro-represents the functor
$\text{\rm Def}_S(\cdot)$ is unique up to unique isomorphism. It is called the universal deformation ring of this functor.
\definition{Proof} For every $A\in\text{Ob}(\Cal C)$ let  $\text{Def}'_S(A)\subseteq\text{Def}(O_l,A)$ denote the subset of isomorphism classes of those deformations $O=((M,\rho),L)$ of $O_l$ which satisfy the condition {\bf D1} above.
Then $\text{Def}_S(A)\subseteq\text{Def}'_S(A)$ and $\text{Def}_S'(\cdot)$ is a sub-functor of $\text{Def}(O_l,\cdot)$.
\enddefinition
\proclaim{Lemma 6.7} It is sufficient to prove that the functor $\text{\rm Def}_S'
(\cdot)$ is  pro-representable by a noetherian local $\Cal O$-algebra.
\endproclaim
\definition{Proof} The demonstration of this lemma uses the same idea as the proof of Proposition 6.1 of [2] on page 323, so we include it only for the reader's convenience. Assume that $\text{Def}_S'(\cdot)$ is pro-representable by a noetherian local $\Cal O$-algebra $R$. Then we are going to show that $\text{Def}_S(\cdot)$ is pro-representable by a quotient $T$ of $R$. For every ideal $\goth a\triangleleft A$ let $\pi_{\goth a}:A\rightarrow A/\goth a$ denote the quotient map. Note that
for every $A\in\text{Ob}(\Cal C)$ and for every $D\in\text{Def}'_S(A)$ the subset $\text{Def}_S(A)\subseteq\text{Def}'_S(A)$ satisfies the following properties:
\roster
\item"$(i)$" we have $f_*(D)\in\text{Def}_S(B)$ for every surjective morphism $f:A\rightarrow B$ in $\Cal C$,
\item"$(ii)$" if $\goth a$ and $\goth b$ are proper ideals of $A$ such that $\pi_{\goth a*}(D)
\in\text{Def}_S(A/\goth a)$ and $\pi_{\goth b*}(D)\in\text{Def}_S(A/\goth bA)$ then $\pi_{\goth a\cap\goth b*}(D)\in\text{Def}_S(A/(\goth a\cap\goth b))$,
\item"$(iii)$" if $f:A\rightarrow B$ is an injective morphism in $\Cal C$ then
we have $D\in\text{Def}_S(A)$ if and only if $f_*(D)\in\text{Def}_S(B)$.
\endroster
Let $J$ denote the set of those open and closed ideals $\goth n$ of $R$ such that the isomorphism class deformations of $O_l$ over $R/\goth n$ corresponding to the quotient map $R\rightarrow R/\goth n$ lies in $\text{Def}_S(R/\goth n)$. Then $J$ is closed under finite intersections by property $(ii)$ and we have $\goth r\in J$ for every open and closed ideal $\goth r$ which contains an element of $J$ by property $(i)$. Let $\goth j$ denote the intersection of all ideals belonging to the set $J$: it is a closed ideal. We claim that every open and closed ideal $\goth n$ which contains $\goth j$ is an element of $J$. Because the complement of $\goth n$ is a compact set it can be covered by a finite collection of open sets which are complements of ideals belonging to $J$ in this case. Therefore $\goth n$ contains the intersection of finitely many elements of $J$ so it must belong to this set by the above. Hence every $\Cal O$-algebra homomorphism from $R$ into an object of $\Cal C$ whose kernel contains $\goth j$ corresponds to an element of $\text{Def}_S(A)$. On the other hand let $A$ be an object of $\text{Ob}(\Cal C)$, let $D$ be an element of $\text{Def}_S(A)$ and let $f:R\rightarrow A$ the $\Cal O$-algebra homomorphism corresponding to $D$. Then $f$ is the composition of the quotient map $\pi:R\rightarrow R/\goth n$ for some open and closed ideal $\goth n$ and an injective homomorphism $R/\goth n\rightarrow A$. By property $(iii)$ above the element of $\text{Def}'_S(R/\goth n)$ corresponding to $\pi$ actually lies in $\text{Def}_S(R/\goth n)$. In particular the kernel $\goth n$ of the homomorphism $f$ contains $\goth j$. Hence the ring $T=R/\goth j$ pro-represents the functor $\text{Def}_S(\cdot)$.\ $\square$
\enddefinition
Let us return to the proof of Theorem 6.6. We are going to use the notation of Definition 6.3. Let $G_S$ denote the Galois group of the maximal separable extension of $F$ in $\overline F$ which is unramified at every prime $\goth q$ not in $S$. For $i=0,1,2,3$ let $G_i$ and $E_i$ denote the kernel of the canonical surjection $p_i:GL_2(A_i)\rightarrow GL_2(\Bbb F_l)$ and the set
$$E_i=\{\rho:G_S\rightarrow GL_2(A_i)|p_i\circ\rho=\rho_l\},$$
respectively. Let moreover $H_i$ denote the set of free $A$-submodules $L\subset A^2$ of rank one which are direct summands of $A^2$ and $L\otimes_A\Bbb F_l=\overline L$. The group $G_i$ acts on $E_i$ by conjugation and on $H_i$ via its regular left action.
Then
$$\text{Def}_S'(A_i)=(E_i\times H_i)/G_i\quad(i=0,1,2,3)$$
where we let $G_i$ act diagonally on the product $E_i\times H_i$. For every pair $\rho\times L\in E_i\times H_i$ let $G_i(\rho\times L)$ denote the stabilizer of $\rho\times L$ in $G_i$. For any $\Cal C$-morphism $f:A_i\rightarrow A_j$ for some $i$ and $j$ let $D(f)$ denote the map $D(f):E_i\times H_i\rightarrow E_j\times H_j$ given by the rule $D(f)(\rho\times L)=
f\circ\rho\times L\otimes_{f(A_i)}A_j$. Then the map $h$ of Definition 6.3 is just the map
$$h:(E_i\times H_i)/G_i\rightarrow(E_i\times H_i)/G_i\times_{(E_i\times L_i)/G_i}
(E_i\times H_i)/G_i$$
induced by $D(\pi_1)$ and $D(\pi_2)$.
\proclaim{Lemma 6.8} Assume that $\phi_1$ is surjective. The following holds:
\roster
\item"$(i)$" the map $h$ is surjective,
\item"$(ii)$" the map $h$ is injective if the homomorphism $G_1(\rho\times L)\rightarrow G_0(D(\phi_1)(\rho\times L))$ induced by $\phi_1$ is surjective for every $\rho\times L\in E_1\times H_1$.
\endroster
\endproclaim
\definition{Proof} Let $\rho_1\times L_1\in E_1\times H_1$ and $\rho_2\times L_2\in E_2\times H_2$ be such that $D(\phi_1)(\rho_1\times L_1)$ and $D(\phi_2)(\rho_2\times L_2)$ lie in the same $G_0$-orbit. The homomorphism $G_1\rightarrow G_0$ induced by $\phi_1$ is surjective hence there is a pair $\rho'_1\times L'_1\in E_1\times H_1$ in the $G_1$-orbit of $\rho_1\times L_1$ such that $D(\phi_1)(\rho'_1\times L'_1)$ is equal to $D(\phi_2)(\rho_2\times L_2)$. Then there is a pair $\rho_3\times L_3\in E_3\times H_3$ such that $D(\pi_1)(\rho_3\times L_3)=\rho_1'\times L_1'$ and $D(\pi_2)(\rho_3\times L_3)=\rho_2\times L_3$ hence claim $(i)$ is true. Now assume that the condition in claim $(ii)$ holds and let $\rho_3'\times L_3'\in E_3\times H_3$ be another pair such that $D(\pi_1)(\rho_3'\times L_3')$ lies in the $G_1$-orbit of
$\rho_1'\times L_1'$ and $D(\pi_2)(\rho_3'\times L_3')$ lies in the $G_2$-orbit of $\rho_2\times L_2$. In order to prove that claim $(ii)$ is true it will be enough to show that $\rho_3\times L_3$ and $\rho_3'\times L_3'$ are in the same $G_3$-orbit. Let $g_1\in G_1$ and $g_2\in G_2$ be two elements such that $D(\pi_1)(\rho_3'\times L_3')g_1=\rho_1'\times L'_1$ and $D(\pi_2)(\rho_3'\times L_3')g_2=\rho_2\times L_2$. Let $h_1$ and $h_2$ denote the image of $g_1$ under the natural map $G_1\rightarrow G_0$ and the image of $g_2$ under the natural map $G_2\rightarrow G_0$, respectively. Then $h_2h_1^{-1}\in G_0(D(\phi_1)(\rho'_1\times L'_1))$ hence there is a $g_1'\in G_1(\rho_1'\times L_1')$ such that the image of $g_1'$ maps to $h_2h_1^{-1}$ under the natural map $G_1\rightarrow G_0$. The pair $g_1'g_1\times g_2$ lies in the fiber product $G_3=G_1\times_{G_0}G_2$ and maps $\rho_3'\times L_3'$ to
$\rho_3\times L_3$.\ $\square$
\enddefinition
With the aid of Lemma 6.8 the proof of Theorem 6.6 is now easy. Condition {\bf H1} is immediate from claim $(i)$ above. On the other hand conditions {\bf H2} and {\bf H4} follow at once from claim $(ii)$ above and Lemma 6.9 below. We only need to prove now that $\dim_{\Bbb F_l}(t_{\text{Def}'_S(\cdot)})$ is finite which is equivalent to the set $\text{Def}'_S(\Bbb F_l[\epsilon])$ being finite. The set of affine lines $L\subset\Bbb F_l[\epsilon]^2$ such that $L\otimes_{\Bbb F_l[\epsilon]}\Bbb F_l=\overline L$ is obviously finite. On the hand the set
$$\{\rho:G_S\rightarrow GL_2(\Bbb F_l[\epsilon])|a\circ\rho=\rho_l\}$$
is well-known to be bijective to the cohomology group $H^1(G_S,\text{Ad}(\rho_l))$ which is finite.\ $\square$
\proclaim{Lemma 6.9} For every pair $\rho\times L\in E_i\times H_i$ the stabilizer $G_i(\rho\times L)$ is equal to $G_i\cap Z(A_i)$ where $Z(A_i)$ denotes the subgroup of scalar matrices in $GL_2(A_i)$.
\endproclaim
\definition{Proof} Let $n$ denote the smallest natural number such that $\goth m_{A_i}^n\neq 0$ but $\goth m_{A_i}^{n+1}=0$. We are going to show that the stabilizer of $\rho\times L$ in $GL_2(A_i)$ is $Z(A_i)$ by induction on $n$. When $n=0$ then $A_i=\Bbb F_l$ and the claim above follows from a straightforward computation. Assume now that we know the claim for $n-1$ and let $r$ be an element of the stabilizer of $\rho\times L$. By induction the image of $r$ under the natural projection $GL_2(A_i)\rightarrow GL_2(A_i/\goth m_{A_i}^n)$ is a scalar matrix. Hence we may assume that $r$ is actually the identity matrix under the image of this projection multiplying by a suitable scalar matrix. Now we only have to prove that $r$ acts on $(m_{A_i}^n)^2\subset A_i^2$ as a scalar multiplication by an element of $1+\goth m_{A_i}^n$. But this is clear from the case $n=0$ as $\rho|_{(m_{A_i}^n)^2}=\rho_l\otimes_{\Bbb F_l}(m_{A_i}^n)^2$ and $L\cap(m_{A_i}^n)^2=\overline L\otimes_{\Bbb F_l}(m_{A_i}^n)^2$.\ $\square$
\enddefinition

\heading 7. Homomorphisms between deformation rings and Hecke rings
\endheading

\definition{Notation 7.1} Let $V=(M,\rho)$ be a representation of the topological group $G$ over an object $A$ of the category $\Cal C$ of Definition 6.2. For every $\Cal C$-morphism $f:A\rightarrow A'$ let $V\otimes_AA'$ denote the couple $(M\otimes_AA',\rho_{A'})$ where $\rho_{A'}$ is the composition of $\rho$ and the map $GL_A(M)\rightarrow GL_{A'}(M\otimes_AA')$ which assigns $\phi\otimes_A\text{id}_{A'}$ to $\phi\in GL_A(M)$. Recall that two representations $(M,\rho)$ and $(M',\rho')$ of $G$ are said to be isomorphic if $M$ and $M'$ are isomorphic as $A[G]$-modules. When $A$ is a field one may attach to every representation $V=(M,\rho)$ another representation $V_{ss}=(M_{ss},\rho_{ss})$, unique up to isomorphism, such that $M_{ss}$ is semi-simple as an $A[G]$-module and has the same Jordan-H\"older components as the $A[G]$-module $M$. Let $\bold k$ denote now a finite extension of $\Bbb F_l$ and let $V=(M,\rho)$ be a representation of $\text{\rm Gal}(\overline F|F)$ over $\bold k$ which is unramified away from $\goth p$ and $\infty$, contains an affine $\bold k$-line $L$ which is fixed by $I_{\goth p}$ but not left stable by $\text{\rm Gal}(\overline F|F)$ and the semi-simplification of $V$ is isomorphic to the semi-simplification of $(\Bbb F_l^2,\rho_l)\otimes_{\Bbb F_l}\bold k$. Suppose moreover that the $\bold k[D_{\infty}]$-module $M$ has a filtration:
$$\CD0@>>>\bold k(1)@>>>M@>>>\bold k@>>>0\endCD$$
where $\bold k(1)$, $\bold k$ denotes the free $\bold k$-module of rank one equipped with the $D_{\infty}$-action given by $\chi_{\bold k}|_{D_{\infty}}$ and the trivial character, respectively.
\enddefinition
\proclaim{Lemma 7.2}  Under the assumptions above we have $V\cong(\Bbb F_l,\rho_l) \otimes_{\Bbb F_l}\bold k$.
\endproclaim
\definition{Proof} Because the Jordan-H\"older components of the $\bold k[\text{Gal}(\overline F|F)]$-modules corresponding to $V$ and $(\Bbb F_l,\rho_l)\otimes_{\Bbb F_l}\bold k$ are the same the $\bold k[\text{Gal}(\overline F|F)]$-module $M$ is the extension of one of the $\bold k[\text{Gal}(\overline F|F)]$-modules given by the representations $(\Bbb F_l,\overline{\chi}_l)\otimes_{\Bbb F_l}\bold k$ or $(\Bbb F_l,1)\otimes_{\Bbb F_l}\bold k$ by the other where $1$ denotes the trivial representation. Because both $\overline{\chi}_l$ and $1$ are unramified at $\goth p$ the $\bold k[\text{Gal}(\overline F|F)]$-module $M$ contains a non-zero submodule which is fixed by the inertia group $I_{\goth p}$. On the other hand it also contains the affine line $L$ which is fixed by $I_{\goth p}$ but not left stable by $\text{Gal}(\overline F|F)$ hence as a $\bold k$-vectorspace $M$ is spanned by vectors which is fixed by $I_{\goth p}$. Therefore $\rho$ is unramified at $\goth p$. In a suitable basis the inertia group $I_{\infty}$ acts by upper-triangular matrices which are equal to $1$ on the diagonal. Therefore the representation $\rho$ is tamely ramified at $\infty$. Because every \'etale cover of the affine line $\Bbb A^1_{\Bbb F_q}$ which is tamely ramified at $\infty$ is in fact unramified at $\infty$ as well the Galois representation $\rho$ is everywhere unramified. The only everywhere unramified extensions of $F$ are the constant field extensions. As a representation of the absolute Galois group $\text{Gal}
(\overline{\Bbb F}_q|\Bbb F_q)$, where the latter group is isomorphic to the profinite completion of $\Bbb Z$, the representation $(\Bbb F_l^2,\rho_l)\otimes_{\Bbb F_l}\bold k$ is the only one up to isomorphism which has an affine $\bold k$-line not fixed by the Galois action and has Jordan-H\"older components $(\Bbb F_l,\overline{\chi}_l)\otimes_{\Bbb F_l}\bold k$ and $(\Bbb F_l,1)\otimes_{\Bbb F_l}\bold k$.\ $\square$
\enddefinition
Let $\text{Def}_{min}(\cdot):\Cal C\rightarrow Sets$ denote the functor $\text{Def}_S(\cdot)$ when $S$ is the empty set. The functor $\text{Def}_{min}(\cdot)$ is called the minimal deformation problem.
\proclaim{Proposition 7.3} The universal deformation ring of the minimal deformation problem $\text{\rm Def}_{min}(\cdot)$ is $\Bbb Z_l$.
\endproclaim
\definition{Proof} We have to show that for every $A\in\text{Ob}(\Cal C)$ the set $\text{Def}_{min}(A)$ has exactly one element because we have exactly one
$\Bbb Z_l$-algebra homomorphism from $\Bbb Z_l$ into $A$. First we are going to show that $\text{Def}_{min}(A)$ is not empty. Let $\text{Fr}_{\infty}$ be again a lift of the Frobenius at the place $\infty$ in $\text{Gal}(\overline F|F)$. Let
$\phi_l:\text{Gal}(\overline F|F)\rightarrow\Bbb Z_l$ be the unique map which factors through the Galois group of the maximal everywhere unramified extension of $F$ such that $\phi_l(\text{Fr}_{\infty})=1$ and $\phi_l(gh)=\chi_l(h)\phi_l(g)+\phi_l(h)$ for every $g$, $h\in\text{Gal}(\overline F|F)$. Recall that $\chi_A$ denotes the composition of the character $\chi_l:\text{Gal}(\overline F|F)\rightarrow\Bbb Z_l^*$ and the natural map $\Bbb Z_l^*\rightarrow A^*$. Moreover let $\phi_A$ be the composition of the twisted homomorphism $\phi_l:\text{Gal}(\overline F|F)\rightarrow\Bbb Z_l$ and the natural map $\Bbb Z_l\rightarrow A$. Let $\rho_A:\text{Gal}(\overline F|F)\rightarrow GL_2(A)$ denote the Galois representation such that
$$\rho_A(g)=\left(\matrix1&0\\0&\chi_A(g)\endmatrix\right)
\quad(\forall g\in\text{Gal}(\overline F|F))$$
when $l$ does not divide $q-1$ and
$$\rho_A(g)=\left(\matrix1&\phi_A(g)\\0&\chi_A(g)\endmatrix\right)
\quad(\forall g\in\text{Gal}(\overline F|F))$$
when $l$ does divide $q-1$. Then for every $A$-module $L\subset A^2$ of rank one such that $L\otimes_A\Bbb F_l=\overline L$ the isomorphism class of the pair $((A^2,\rho_A),L)\in\Cal P_{\text{Gal}(\overline F|F)}(A)$ lies in $\text{Def}_{min}(A)$ as $\rho_A$ is unramified at $\goth p$. Let $O=((M,\rho),L)$ be now an arbitrary element of $\text{Def}_{min}(A)$. By condition {\bf D4} the inertia group $I_{\infty}$ acts by upper-triangular matrices which are equal to $1$ on the diagonal in a suitable $A$-basis of $M$. As the additive group of $A$ is $l$-primary group the representation $\rho$ is tamely ramified at $\infty$. Hence it must be everywhere unramified by the argument presented in the proof of Lemma 7.2 above. In particular the image of $\rho$ is generated by $\rho(\text{Fr}_{\infty})$. In any $A$-basis furnished by condition {\bf D4} the matrix of $\rho(\text{Fr}_{\infty})$ is of the form:
$$\rho(\text{Fr}_{\infty})=\left(\matrix1&x\\0&\chi_A(g)\endmatrix\right)$$
for some $x\in A$. Because the vector $(0,1)$ is fixed by $\chi_A^{-1}\rho(\text{Fr}_{\infty})$, its reduction modulo $\goth m_A$ spans a one-dimensional subspace of $M\otimes_A\Bbb F_l$ fixed by $\chi^{-1}_A\rho$ mod $\goth m_A$. But such a subspace is unique as $\rho$ mod $\goth m_A$ is a conjugate of $\overline{\rho}_l$. Hence for a suitable $y\in A^*$ the matrix of $\rho(\text{Fr}_{\infty})$ mod $\goth m_A$ in the basis $(y,0),(0,1)$ is the same as the matrix of $\rho_A(\text{Fr}_{\infty})$. Therefore there is a matrix $m\in GL_2(A)$ such that $m$ is congruent to the identity matrix mod $\goth m_A$ and $m^{-1}\rho(\text{Fr}_{\infty})m=\rho_A(\text{Fr}_{\infty})$. In other words we may assume that $(M,\rho)=(A^2,\rho_A)$ and $L\otimes_A\Bbb F_l=\overline L$. Let $H\subset GL_2(A)$ denote the subgroup:
$$H=\{\left(\matrix 1&0\\0&y\endmatrix\right)\in GL_2(A)|y\equiv1\!\!\!\mod\goth m_A\}$$
when $l$ does not divide $q-1$ and let $H$ denote the subgroup:
$$H=\{\left(\matrix 1&x\\0&1+(q-1)x\endmatrix\right)\in GL_2(A)|x\equiv0\!\!\!\mod\goth m_A\}$$
when $l$ does divide $q-1$. Then $H$ commutes with $\rho_A(\text{Fr}_{\infty})$ and acts transitively on the set of affine lines $L\subset A^2$ such that $L\otimes_A\Bbb F_l=\overline L$. Hence the set $\text{Def}_{min}(A)$ has at most one
element.\ $\square$
\enddefinition
\definition{Definition 7.4} For every commutative ring $A$ with unity and for every $A$-linear endomorphism $m$ of a finitely generated free $A$-module let $\text{Tr}(m)\in A$ denote the trace of $m$. Let $\text{Def}_{\goth p}(\cdot):\Cal C\rightarrow Sets$ denote the functor $\text{Def}_S(\cdot)$ when $S$ is the one-element set $\{\goth p\}$ and let $R(\goth p)$ denote its universal deformation ring. The minimal deformation functor $\text{Def}_{min}(\cdot)$ is a sub-functor of $\text{Def}_{\goth p}(\cdot)$. Let $\pi_R:R(\goth p)\rightarrow\Bbb Z_l$ be the surjective $\Bbb Z_l$-algebra homomorphism of universal deformation rings corresponding to this inclusion and let $I_R\triangleleft R(\goth p)$ be the kernel of this projection. For every place $v$ of $F$ let $\text{Fr}_v$ be a lift of the Frobenius at the place $v$ in $\text{Gal}(\overline F|F)$. Let $\goth n\triangleleft R(\goth p)$ be an open and closed ideal and let $((M,\rho),L)\in\Cal P_{\text{Gal}(\overline F|F)}
(R(\goth p)/\goth n)$ be a pair whose isomorphism class corresponds to the element of $\text{Def}_{\goth p}(R(\goth p)/\goth n)$ given by the canonical projection $R(\goth p)\rightarrow R(\goth p)/\goth n$. Then for every prime ideal $\goth q\neq\goth p$ of $\bold A$ the trace $\text{Tr}(\rho(\text{Fr}_{\goth q}))
\in R(\goth p)/\goth n$ is independent of choice of $\text{Fr}_{\goth q}$ as $\rho$ is unramified at $\goth q$ (and of course it depends only on the isomorphism class of $((M,\rho),L)$). The elements $\text{Tr}(\rho(\text{Fr}_{\goth q}))$ as $\goth n$ runs through the set of open and closed ideals of $R(\goth p)$ form a compatible family and therefore has a limit in $R(\goth p)$ which will be denoted by $\text{\rm Tr}(\rho_{univ}(\text{\rm Fr}_{\goth q}))$.
\enddefinition
\proclaim{Proposition 7.5} The ideal $I_R$ is generated by the set:
$$\{1+q^{\deg(\goth q)}-\text{\rm Tr}(\rho_{univ}(\text{\rm Fr}_{\goth q}))|\goth p\neq
\goth q\triangleleft\bold A\text{\ is a prime ideal.}\}$$
\endproclaim
\definition{Proof} (Compare with Proposition 3.3 of [1], pages 109-111.) Let $J$ denote the ideal generated by the set in the claim above. Clearly $J$ is a subset of $I_R$ so we only have to show the other inclusion. Because $R(\goth p)$ is a noetherian complete local ring the ideal $J$ is the intersection of all open and closed ideals which contain $J$. Therefore it will be sufficient to prove that for every object $A$ of $\Cal C$ and for every pair $O=((M,\rho),L)$ deforming $O_l$ whose isomorphism class lies in $\text{Def}_{\goth p}(A)$ such that the image of $J$ with respect to the $\Bbb Z_l$-homomorphism $R(\goth p)\rightarrow A$ corresponding to $O$ is zero the Galois representation $\rho$ is unramified at $\goth p$. Assume first that $l$ does not divide $q-1$. By condition {\bf D4} the
module $M$ has an $A$-basis such that the matrix of $\rho(\text{Fr}_{\infty})$ is upper-triangular and the diagonal element of the first and second row is $1$ and $q$, respectively. Conjugating by a suitable upper-triangular matrix we may even assume that $\rho(\text{Fr}_{\infty})$ is of the form:
$$\rho(\text{Fr}_{\infty})=\left(\matrix1&0\\0&q\endmatrix\right).$$
If we write
$$\rho(g)=\left(\matrix a(g)&b(g)\\c(g)&d(g)\endmatrix\right)$$
for every $g\in\text{Gal}(\overline F|F)$ in this basis then
$$a(g)={1\over q-1}(q\text{Tr}(\rho(g))-\text{Tr}(\rho(\text{Fr}_{\infty}g)))$$
and
$$d(g)={1\over q-1}(\text{Tr}(\rho(\text{Fr}_{\infty}g))-\text{Tr}(\rho(g))).$$
By assumption $\text{Tr}(\text{Fr}_{\goth q})=1+\chi_A(\text{Fr}_{\goth q})$ for every  prime ideal $\goth q\neq\goth p$ of $\bold A$ hence $\text{Tr}(h)=1+\chi_A(h)$ for every $h\in\text{Gal}(\overline F|F)$ as well by the Chebotarev density theorem. Therefore the equations above imply that $a(g)=1$ and $d(g)=\chi_A(g)$ because $\chi_A(\text{Fr}_{\infty})=q$. In particular for every $g\in I_{\goth p}$ we have:
$$\rho(g)=\left(\matrix 1&b(g)\\c(g)&1\endmatrix\right).$$
The affine line spanned by the vectors $(1,0)$ and $(0,1)$ modulo $\goth m_A$ are fixed by $\rho(\text{Fr}_{\infty})$ mod $\goth m_A$ hence by the action of the whole absolute Galois group as well. Hence the free $A$-module $L\subset A^2$ of rank one is generated by a vector of the form $(1,x)$ for some $x\in A^*$. Because
$$\rho(g)(1,x)=(1+b(g)x,c(g)+x)$$
for every $g\in I_{\goth p}$ we must have $b(g)=c(g)=0$ in this case. In particular $\rho$ is unramified at $\goth p$. Assume now that $l$ divides $q-1$. Using condition {\bf D4} as above we may conclude that the module $M$ has an $A$-basis such that the matrix of $\rho(\text{Fr}_{\infty})$ is of the form:
$$\rho(\text{Fr}_{\infty})=\left(\matrix1&1\\0&q\endmatrix\right)$$
in this basis. We may even assume that $L$ is spanned by the vector $(0,1)$. If we write
$$\rho(g)=\left(\matrix a(g)&b(g)\\c(g)&d(g)\endmatrix\right)$$
for every $g\in\text{Gal}(\overline F|F)$ in this basis then
$$\text{Tr}(\rho(\text{Fr}_{\infty}g))=a(g)+c(g)+qd(g)=1+q\chi_A(g)$$
using the same reasoning as above. Since every $g\in I_{\goth p}$ fixes $(0,1)$ we get that $b(g)=0$ and $d(g)=1$. The determinant of $\rho(g)$ is equal to $\chi_A(g)=1$ hence $a(g)=1$. Therefore $c(g)=1+q\chi_A(g)-(1+q)=0$ by the equation above. So $\rho$ is unramified at $\goth p$ in this case as well.\ $\square$
\enddefinition
The following corollary follows from Proposition 7.5 above the same way as Corollary 3.4 of [1] on page 111 from Proposition 3.3 of the same paper.
\proclaim{Corollary 7.6} The complete local $\Bbb Z_l$-algebra $R(\goth p)$ is topologically generated by the elements $\text{\rm Tr}(\rho_{univ}(\text{\rm Fr}_{\goth q}))$ where $\goth q\neq\goth p$ is any proper prime ideal of $\bold A$.\ $\square$
\endproclaim
\definition{Definition 7.7} Let $\widehat{\Bbb T}(\goth p)$ denote the commutative
$\Bbb Z$-algebra with unity generated by the endomorphisms $T_{\goth q}$ of
the $\Bbb Z$-module $H(\Cal T,\Bbb Z)^{\Gamma_0(\goth p)}$ where $\goth q\triangleleft \bold A$ is any prime ideal different from $\goth p$. By Corollary 3.13 the algebra $\widehat{\Bbb T}(\goth p)$ is a sub-algebra of the endomorphism ring of a finitely generated, free $\Bbb Z$-module hence it must be a finitely generated, free $\Bbb Z$-module, too. Let $\widehat{\goth E}(\goth p)$ denote the ideal of $\widehat{\Bbb T}(\goth p)$ generated by the elements $T_{\goth q}-q^{\deg(\goth q)}-1$, where $\goth
q\neq\goth p$ is any prime. Because every generator $T_{\goth q}$ of $\Bbb
T_l(\goth p)$ is congruent to an element of $\Bbb Z$ modulo the ideal $\widehat{\goth E}(\goth p)$ the natural inclusion of $\Bbb Z$ in $\widehat{\Bbb T}(\goth p)$ induces a surjection $\Bbb Z@>>>\widehat{\Bbb T}(\goth p)/\widehat{\goth E}(\goth p)$. This map is also injective as the annihilator of the action of $\widehat{\Bbb T}(\goth p)$ on the subgroup of $H(\Cal T,\Bbb Z)^{\Gamma_0(\goth p)}$ generated by $E_{\goth p}$ is equal to $\widehat{\goth E}(\goth p)$. Therefore $(l,\widehat{\goth E}(\goth p))$ is a maximal ideal in $\widehat{\Bbb T}(\goth p)$ with residue field $\Bbb F_l$. Let $T(\goth p)$ denote the completion of $\widehat{\Bbb T}(\goth p)$ with respect to the maximal ideal $(l,\widehat{\goth E}(\goth p))$ and for every prime ideal $\goth q\neq\goth p$ of $\bold A$ let $T_{\goth q}\in T(\goth p)$ denote the image of the Hecke operator $T_{\goth q}$ under the canonical projection $\widehat{\Bbb T}(\goth p)\rightarrow T(\goth p)$ by slight abuse of notation. The ring $T(\goth p)$ is naturally equipped with the structure of a local $\Bbb Z_l$-algebra.
\enddefinition
\proclaim{Proposition 7.8} There is a unique $\Bbb Z_l$-algebra homomorphism $\phi:R(\goth p)\rightarrow T(\goth p)$ such that
$$\phi(\text{\rm Tr}(\rho_{univ}(\text{\rm Fr}_{\goth q})))=T_{\goth q}$$
for every proper prime ideal $\goth q\neq\goth p$ of $\bold A$.
\endproclaim
\definition{Proof} Note that the map $\phi$, if it exits, is automatically surjective and local. So it is unique because it is determined on a set of topological generators by Corollary 7.6. Let $\widetilde T(\goth p)$ denote the normalization of the algebra $T(\goth p)$ and let $\iota:T(\goth p)\rightarrow\widetilde T(\goth p)$ be the normalization map. Because the $T(\goth p)$ is a finitely generated free $\Bbb Z_l$-module the $\Bbb Z_l$-algebra $\widetilde T(\goth p)$ can be written as a finite product $\prod_{i\in I}\Cal O_i$ where each $\Cal O_i$ is a discrete valuation ring which is finitely generated and free as a $\Bbb Z_l$-module. For every $i\in I$ let $\pi_i:T(\goth p)\rightarrow\Cal O_i$ be the composition of $\iota$ and the canonical projection. For every $i\in I$ we are going to construct a
$\Bbb Z_l$-algebra homomorphism $\phi_i:R(\goth p)\rightarrow\Cal O_i$ such that
$$\phi_i(\text{\rm Tr}(\rho_{univ}(\text{\rm Fr}_{\goth q})))=\pi_i(T_{\goth q})\tag"(7.8.1)"$$
for every proper prime ideal $\goth q\neq\goth p$ of $\bold A$. In this case the product
$$\phi=\prod_{i\in I}\phi_i:R(\goth p)\rightarrow\prod_{i\in I}\Cal O_i=
\widetilde T(\goth p)$$
maps $R(\goth p)$ into the image of $T(\goth p)$ with respect to $\iota$ because $\phi$ maps
$\text{\rm Tr}(\rho_{univ}(\text{\rm Fr}_{\goth q})))$ to $\iota(T_{\goth q})$ for every proper prime ideal $\goth q\neq\goth p$ of $\bold A$. But $\iota$ is injective hence the claim above implies the proposition.

Let $\goth m_i$, $\bold k_i$ denote the maximal ideal and the residue field of $\Cal O_i$, respectively. Let $\Cal O'_i$ denote the preimage of $\Bbb F_l\subseteq\bold k_i$ with respect to the residue map mod $\goth m_i$. Then $\Cal O'_i/\goth m_i=\Bbb F_l$ hence $\Cal O'_i/\goth m_i^n$ is an object of the category $\Cal C$ for every positive integer $n$. Let $\pi_{in}:GL_2(\Cal O_i)\rightarrow GL_2(\Cal O_i/\goth m_i^n)$ denote the canonical projection. Note that it is sufficient to prove that there is a representation $\rho_i:\text{Gal}(\overline F|F)\rightarrow GL_2(\Cal O_i)$ continuous with respect to the $\goth m_i$-adic topology on $GL_2(\Cal O_i)$ and the Krull topology on $\text{Gal}(\overline F|F)$ along with a free $\Cal O_i$-module $L_i\subset\Cal O_i^2$ of rank one satisfying the following properties:
\roster
\item"{\bf L1.}" the order pair $((\Bbb F_l^2,\pi_{i1}\circ\rho_i),L_i\otimes_{\Cal O_i}\Bbb F_l)$ is isomorphic to $O_l\otimes_{\Bbb F_l}\bold k_i$,
\item"{\bf L2.}" representation $\rho_i$ is unramified at every prime ideal $\goth q\neq\goth p$ of $\bold A$,
\item"{\bf L3.}" the inertia subgroup $I_{\goth p}$ at $\goth p$ acts trivially on the submodule $L_i$,
\item"{\bf L4.}" the determinant of $\rho_i$ is equal to the composition $\chi_i$ of the cyclotomic character $\chi_l:\text{Gal}(\overline F|F)\rightarrow\Bbb Z_l^*$ and the natural map $\Bbb Z_l^*\rightarrow\Cal O_i^*$,
\item"{\bf L5.}" in a suitable $\Cal O_i$-basis the matrix of $\rho_i(g)$ for every $g\in D_{\infty}$ is of the form:
$$\rho_i(g)=\left(\matrix1&\psi_i(g)\\0&\chi_i(g)\endmatrix\right)$$
for some $\psi_i(g)\in\Cal O_i$,
\item"{\bf L6.}" we have $\text{Tr}(\rho_i(\text{Fr}_{\goth q}))=\pi_i(T_{\goth q})$ for every prime ideal $\goth q\neq\goth p$ of $\bold A$.
\endroster
In fact when the claim above holds we may assume that the image of $\rho_i$ lies in $GL_2(\Cal O'_i)$ and $L_i=L_i'\otimes_{\Cal O'_i}\Cal O_i$ for some free $\Cal O_i'$-module $L_i'\subset
(\Cal O_i')^2$ of rank one without loss of generality. In this case the isomorphism class of the pair $(((\Cal O'_i/\goth m_i)^2,\pi_{in}\circ\rho_i),L'_i\otimes_{\Cal O'_i}\Cal O'_i/\goth m_i^n)$ lies in $\text{Def}_{\goth p}(\Cal O'_i/\goth m_i^n)$ for every positive integer $\goth n$. Let $\phi_{in}:R(\goth p)\rightarrow\Cal O'_i/\goth m_i^n$ denote the corresponding map. The maps $\phi_{in}$ satisfy the obvious compatibility and their limit $\phi:R(\goth p)\rightarrow \Cal O'_i$ satisfies the equation (7.8.1) above by condition {\bf L6}.

In the rest of the proof we occupy ourselves with the proof of the existence of $\rho_i$ and $L_i$. Assume first that the kernel of the projection $\pi_i:T(\goth p)
\rightarrow\Cal O_i$ contains the image of the ideal $\widehat{\goth E}(\goth p)$ with respect to the canonical projection $\widehat{\Bbb T}(\goth p)\rightarrow T(\goth p)$. Then $\Cal O_i=\Bbb Z_l$ and $\pi_i(T_{\goth q})=1+q^{\deg(\goth q)}$ for every
prime ideal $\goth q\neq\goth p$ of $\bold A$. Clearly the construction of Proposition 7.3 supplies the representation $\rho_i$ in this case. Assume next that the kernel of the projection $\pi_i:T(\goth p)\rightarrow\Cal O_i$ does not contain the image of the ideal $\widehat{\goth E}(\goth p)$ with respect to the canonical projection $\widehat{\Bbb T}(\goth p)\rightarrow T(\goth p)$. Let $j:\widehat{\Bbb T}(\goth p)\rightarrow\Bbb T(\goth p)\oplus\Bbb Z$ denote the direct sum of the surjections $j_1:\widehat{\Bbb T}(\goth p)\rightarrow\Bbb T(\goth p)$ and $j_2:\widehat{\Bbb T}(\goth p)\rightarrow\Bbb Z$ induced by the restriction of the action of $\widehat{\Bbb T}(\goth p)$ onto $H_!(\Cal T,\Bbb Z)^{\Gamma_0(\goth p)}$ and $\Bbb ZE_{\goth p}$, respectively. By Theorem 3.12 and its Corollary 3.13 the direct sum $H_!(\Cal T,\Bbb Z)^{\Gamma_0(\goth p)}\oplus\Bbb ZE_{\goth p}$ is a subgroup of finite index of $H(\Cal T,\Bbb Z)^{\Gamma_0(\goth p)}$. Hence $j\otimes\Bbb Q:\widehat{\Bbb T}(\goth p)\otimes\Bbb Q\rightarrow\Bbb T(\goth p)\otimes\Bbb Q\oplus\Bbb Q$ is an isomorphism. Therefore the composition $\pi_i'$ of the canonical surjection $\widehat{\Bbb T}(\goth p)\rightarrow T(\goth p)$ and $\pi_i$ factors through $j_1$.

The image $R$ of $\pi_i'$ is an order in a number field $K$. The image of the ideal $(\widehat{\goth E}(\goth n),l)$ with respect to $\pi_i'$ is a maximal ideal $\goth m$ in $R$ such that $R/\goth m=\Bbb F_l$. Let $K_i$ denote the completion of $K$ with respect to the valuation corresponding to the ideal $\goth m$. The valuation ring of $K_i$ is obviously $\Cal O_i$. Let $A$ be an abelian variety over $F$ equipped with $R$-multiplication whose existence is guaranteed by Theorem 4.8. The $R$-multiplication on $A$ induces a $K_i$-linear structure on $V_l(A)$ which makes the latter a $K_i$-vectorspace of dimension two. Let $W$ denote this vector space. The action of the absolute Galois group of $F$ on $V_l(A)$ is $K_i$-linear hence it induces a homomorphism $\rho_i':\text{Gal}(\overline F|F)\rightarrow GL_{K_i}(W)$. Because the topological group $\text{Gal}(\overline F|F)$ is compact and the representation is continuous with respect to the topology of $W$ induced by the valuation of $K_i$ there is a free $\Cal O_i$-submodule $V\subset W$ of rank two left stable by $\text{Gal}(\overline F|F)$. Because $A$ has semi-stable reduction at $\goth p$ there is a unique $1$-dimensional $K_i$-subspace $U$ of $W$ fixed by $I_{\goth p}$. Let $L$ denote the intersection of $V$ and $U$ and fix an $\Cal O_i$-linear isomorphism $c:\Cal O_i^2\rightarrow V$. Such an isomorphism furnishes a homomorphism $\rho_i:\text{Gal}(\overline F|F)\rightarrow GL_2(\Cal O_i)$ by restricting the action of $\text{Gal}(\overline F|F)$ via $\rho_i'$ onto $V$. We claim that $\rho_i$ along with $L_i=c^{-1}(L)$ satisfy {\bf L1-L6} for a suitable choice of $V$.

In any case $\rho_i$ satisfies condition {\bf L2} because $A$ has good reduction at every prime ideal $\goth q\neq\goth p$ of $\bold A$. Condition {\bf L3} is obvious and {\bf L4} is consequence of the fact that $\rho_i$ comes from the Tate-module of an abelian variety. Similarly {\bf L5} and {\bf L6} are automatic as $A$ has split multiplicative reduction at $\infty$ and by condition $(iii)$ of Theorem 4.8, respectively. Let $u\in\Cal O_i$ be a uniformizer. Note that
$$\text{Tr}(\rho_i(\text{Fr}_{\goth q}))=\pi_i(T_{\goth q})\equiv1+q^{\deg(\goth q)}\mod\goth m_i$$
for every prime ideal $\goth q\neq\goth p$ of $\bold A$. So the eigenvalues of $\pi_{i1}\circ\rho_i(\text{Fr}_{\goth q})$ are $1$ and $\chi_l(\text{Fr}_{\goth q})$ therefore the semisimplification of $(\bold k_i^2,\pi_{i1}\circ\rho_i)$ is isomorphic to the semisimplification of $(\Bbb F_l^2,\rho_l)\otimes_{\Bbb F_l}\bold k_i$ by the Chebotarev density theorem. Hence the pair $(\rho_i,L_i)$ will also satisfy {\bf L1} if the image of $L$ in $V/uV$ is not left invariant by $\text{Gal}(\overline F|F)$ according to Lemma 7.2. Because $W$ is irreducible as a $\text{Gal}(\overline F|F)$-module the $\Cal O_i$-module $L$ is not left stable by $\text{Gal}(\overline F|F)$. Hence there is a positive integer $n>0$ such that the image of $L$ under the action of $\text{Gal}(\overline F|F)$ does not lie in $u^nV+L$. Let $n$ denote actually the smallest such number and let $V'$ denote the $\Cal O_i$-submodule $u^{n-1}V+L$. Then $V'$ is left invariant by $\text{Gal}(\overline F|F)$ and the image of $L\subset V'$ in $V'/uV'$ is not left invariant by $\text{Gal}(\overline F|F)$.\ $\square$
\enddefinition

\heading 8. Calculations with infinitesimal deformations
\endheading

\definition{Definition 8.1} In this chapter every sheaf or cohomology group is understood to be over the \'etale site. Let $U$ denote the affine curve $\Bbb P^1_{\Bbb F_q}-\{\goth p,\infty\}$ over $\Bbb F_q$. Let $\overline S$ denote the base change of $S$ to $\text{Spec}(\Bbb F_q)$ for every Spec$(\Bbb F_q)$-scheme $S$. Fix a $\overline{\Bbb F}_q$-valued point $v$ of $\overline U$. Then there is an exact sequence:
$$\CD 0@>>>\pi_1^{et}(\overline U,v)@>>>\pi_1^{et}(U,v)@>>>\text{Gal}
(\overline{\Bbb F}_q|\Bbb F_q)@>>>0\endCD$$
of profinite groups. The absolute Galois group $\text{Gal}
(\overline{\Bbb F}_q|\Bbb F_q)$ is isomorphic to the profinite completion $\widehat{\Bbb Z}$ of $\Bbb Z$ and it is equipped with a natural topological generator
$\sigma\in\text{Gal}(\overline{\Bbb F}_q|\Bbb F_q)$, the arithmetic Frobenius. Let $\bold K_l\triangleleft\pi_1^{et}(\overline U,v)$ denote the intersection of the kernels of all surjective homomorphisms $s:\pi_1^{et}(\overline U,v)\rightarrow M$ where $M$ is a finite abelian group of $l$-power order. Then $\bold K_l$ is a characteristic subgroup of $\pi_1^{et}(\overline U,v)$ hence it is a normal subgroup of $\pi_1^{et}(U,v)$. Let $\bold G_l$ denote the quotient $\pi_1^{et}(U,v)/\bold K_l$.
\enddefinition
\proclaim{Lemma 8.2} The following holds:
\roster
\item"$(i)$" the group $\bold G_l$ is isomorphic to a semidirect product
$\bold H_l\rtimes\text{\rm Gal}(\overline{\Bbb F}_q|
\Bbb F_q)$ where $\bold H_l$ is a finitely generated and free $\Bbb Z_l$-module,
\item"$(ii)$" the module $\bold H_l$ has a $\Bbb Z_l$-basis $e_1,e_2,\ldots,e_{\deg(\goth p)}$ such that the action of $\sigma$ on $\bold H_l$ induced by conjugation is given by the rule $\sigma(e_{\deg(\goth p)})=qe_1$ and $\sigma(e_n)=qe_{n+1}$ when $n<\deg(\goth p)$.
\endroster
\endproclaim
\definition{Proof} Let $\bold H_l$ and $t$ denote the image of $\pi_1^{et}(\overline U,v)$ and the image of the Frobenius $\text{Fr}_{\infty}$ in $\bold G_l$ with respect to the quotient map $\pi_1^{et}(U,v)\rightarrow\bold G_l$, respectively. Then $\bold H_l$ is a normal subgroup of $\bold G_l$ and $t\in\bold G_l$ is an element whose image with respect to the canonical surjection $\bold G_l\rightarrow\text{\rm Gal}(\overline{\Bbb F}_q|\Bbb F_q)$ is $\sigma$. Hence the subgroup $\langle t\rangle$ of $\bold G_l$ generated topologically by $t$ is isomorphic to $\text{\rm Gal}(\overline{\Bbb F}_q|\Bbb F_q)$ and $\bold G_l$ is the semidirect product $\bold H_l\rtimes\langle t\rangle$. By definition $\bold H_l$ is the maximal pro-$l$ abelian quotient of $\pi_1^{et}(\overline U,v)$ hence it is a free, finitely generated $\Bbb Z_l$-module isomorphic to $\text{Hom}_{\Bbb Z_l}(H^1(\overline U,\Bbb Z_l),\Bbb Z_l)$. Consider the long cohomological exact sequence:
$$\CD H^0(\overline U,\Cal O^*)@>l^n>>H^0(\overline U,\Cal O^*)@>\delta_n>>
H^1(\overline U,\mu_{l^n})@>>>H^1(\overline U,\Cal O^*)\endCD$$
attached to the Kummer exact sequence:
$$\CD0@>>>\mu_{l^n}@>>>\Cal O^*@>l^n>>\Cal O^*@>>>0\endCD$$
where $\mu_{l^n}$ denotes the sheaf of $l^n$-th roots of unity as usual. The group $H^1(\overline U,\Cal O^*)$ is isomorphic to the Picard group of $\overline U$ but the latter is zero because $\overline U$ is the spectrum of a unique factorization domain. Hence the coboundary maps $\delta_n$ are surjective and their limit induce an isomorphism:
$$\delta:H^1(\overline U,\Cal O^*)\otimes\Bbb Z_l\cong H^1(U,\Bbb Z_l(1)).$$
Let $\Delta$ denote the set of closed points of the projective line over $\overline{\Bbb F}_q$ in the complement of $\overline U$ which are different from $\infty$. Let $\Bbb Z_l[\Delta]$ denote the $\Bbb Z_l$-module generated freely by the elements of $\Delta$. The natural $\text{\rm Gal}(\overline{\Bbb F}_q|\Bbb F_q)$-action on $\Delta$ induces a $\text{\rm Gal}(\overline{\Bbb F}_q|\Bbb F_q)$-module structure on $\Bbb Z_l[\Delta]$. The valuation maps at the elements of $\Delta$ induce a $\text{\rm Gal}(\overline{\Bbb F}_q|\Bbb F_q)$-module isomorphism:
$$j_0:H^0(\overline U,\Cal O^*)\otimes\Bbb Z_l\cong\Bbb Z_l[\Delta].$$
Because $\goth p$ is a prime ideal the set $\Delta$ can be indexed by the natural numbers $n=1,2,\cdots,\deg(\goth p)$ such that $\sigma(n)=n+1$ when $n<\deg(\goth p)$ and $\sigma(\deg(\goth p))=1$. Therefore the group:
$$\bold H_l\cong \text{Hom}_{\Bbb Z_l}(H^1(\overline U,\Bbb Z_l),\Bbb Z_l)=\text{Hom}_{\Bbb Z_l}(H^1(\overline U,\Bbb Z_l(1))\otimes_{\Bbb Z_l}\Bbb Z_l(-1),\Bbb Z_l)$$
satisfies the properties in claims $(i)$ and $(ii)$.\ $\square$
\enddefinition
\definition{Definition 8.3} Any $\Bbb Z_l$-basis $e_1,e_2,\ldots,e_{\deg(\goth p)}$ of $\bold H_l$ where the action of $\sigma$ is given by the rule in part $(ii)$ will be called a $\sigma$-cyclic basis. By slight abuse of notation let $\text{Fr}_x$ denote the image of the Frobenius $\text{Fr}_x$ in $\bold G_l$ with respect to the quotient map $\pi_1^{et}(U,v)\rightarrow\bold G_l$ for every place $x$ of $F$. The image of the inertia group $I_x$ with respect to the quotient map $\pi_1^{et}(U,v)\rightarrow\bold G_l$ lies in $\bold H_l$ hence the induced map $I_x\rightarrow\bold H_l$ factors through the largest pro-$l$ quotient $I_{x,l}$ of $I_x$. The group $I_{x,l}$ is isomorphic to $\Bbb Z_l$ and it has a topological generator $i_x$ such that under the action of $\text{Fr}_x$ on $I_{x,l}$ induced by conjugation the element $i_x$ maps to $q^{\deg(x)}i_x$. By the usual abuse of notation let $i_x$ denote the image of the element $i_x$ in $\bold H_l$ as well with respect to the map $I_{x,l}\rightarrow\bold H_l$ induced by the quotient map $\pi_1^{et}(U,v)\rightarrow\bold G_l$ for every place $x$ of $F$.
\enddefinition
\proclaim{Lemma 8.4} There is a $\sigma$-cyclic basis $e_1,e_2,\ldots,e_{\deg(\goth p)}$ of $\bold H_l$ such that $i_{\goth p}=e_1$ and $i_{\infty}=u\sum_{n=1}^{\deg(\goth p)}e_n$ for some $u\in\Bbb Z_l^*$.
\endproclaim
\definition{Proof} The image of $\text{Fr}_{\goth p}$ with respect to the quotient map $\bold G_l\rightarrow\text{\rm Gal}(\overline{\Bbb F}_q|\Bbb F_q)$ is $\sigma^{\deg(\goth p)}$. Hence $\sigma^{\deg(\goth p)}(i_{\goth p})=q^{\deg(\goth p)}i_{\goth p}$ under the action induced by conjugation. Let $B_l\triangleleft\bold G_l$ denote the closed normal subgroup generated by $i_{\goth p}$. Then $B_l$ is a subgroup of $\bold H_l$ and as a $\Bbb Z_l$-module it is generated by the elements $e_n=q^{1-n}\sigma^{n-1}(i_{\goth p})$ where $n=1,2,\ldots,\deg(\goth p)$. The quotient $\bold G_l/B_l$ is a Galois group which is unramified at $\goth p$ and tamely ramified at $\infty$ hence it is everywhere unramified. Therefore $B_l=\bold H_l$ and $e_1,e_2,\ldots,e_{\deg(\goth p)}$ is a $\sigma$-cyclic basis. The image of $\text{Fr}_{\infty}$ with respect to the quotient map $\bold G_l\rightarrow\text{\rm Gal}(\overline{\Bbb F}_q|\Bbb F_q)$ is $\sigma$ hence $\sigma(i_{\infty})=q
i_{\infty}$ under the action induced by conjugation. Therefore $i_{\infty}$ lies in the $\Bbb Z_l$-module generated by the element $\sum_{n=1}^{\deg(\goth p)}e_n$. The quotient of $\bold H_l$ by the $\Bbb Z_l$-module generated by $i_{\infty}$ is isomorphic to $\pi_1^{et}(\overline S,v)=\text{Hom}_{\Bbb Z_l}(H^1(\overline S,\Bbb Z_l),\Bbb Z_l)$ where $S$ is the complement of $\goth p$ in the projective line. The group $H^1(\overline S,\Bbb Z_l)$ is torsion-free so $i_{\infty}=u\sum_{n=1}^{\deg(\goth p)}e_n$ for some $u\in\Bbb Z_l^*$.\ $\square$
\enddefinition
\definition{Definition 8.5} Let $\text{Def}_{\goth p}^0(A[\epsilon])\subseteq
\text{Def}_{\goth p}(A[\epsilon])$ denote the pre-image of the unique element of $\text{Def}_{min}(A)$  with respect to the map $\text{Def}_{\goth p}(A[\epsilon])
\rightarrow \text{Def}_{\goth p}(A)$ induced by the augmentation homomorphism $a:A[\epsilon]\rightarrow A$. We say that an ordered pair $((A[\epsilon]^2,\rho),L)$ representing an element of $\text{Def}_{\goth p}^0(A[\epsilon])$ is in standard form if the following holds:
\roster
\item"{\bf S1.}" the $A[\epsilon]$-module $L$ is generated by $(1,1)$, when $l$ does not divide $q-1$ and it is generated by $(0,1)$, otherwise,
\item"{\bf S2.}" the matrix of $\rho(g)$ for every $g\in D_{\infty}$ is of the form:
$$\rho(g)=\left(\matrix1&\psi(g)\\0&\chi_{A[\epsilon]}(g)\endmatrix\right)$$
for some $\psi(g)\in A[\epsilon]$ such that $\psi(\text{Fr}_{\infty})=0$, when $l$ does not divide $q-1$ and $\psi(\text{Fr}_{\infty})=1$, otherwise.
\endroster
Because there is no ambiguity in this case we will simply let $\rho$ denote this ordered pair. In order to avoid the use of awkward terminology we will call such an object $\rho$ a representation in standard form. The argument of Proposition 7.5 shows that every element of $\text{Def}_{\goth p}^0(A[\epsilon])$ can be represented by a representation $\rho$ in standard form. Note that the image of $\pi_1^{et}(\overline U,v)$ under $\rho$ is a finite abelian group of $l$-power order. Hence $\rho$ factors through $\bold G_l$. Let $\rho:\bold G_l\rightarrow GL_2(A[\epsilon])$ denote the corresponding homomorphism as well by slight abuse of notation.
\enddefinition
\proclaim{Lemma 8.6} Assume that $\rho$ is in standard form. Then
$$\rho(i_{\goth p})=\cases\left(\matrix 1+a(\rho)\epsilon&-a(\rho)\epsilon\\
a(\rho)\epsilon&1-a(\rho)\epsilon\endmatrix\right),&\text{if $l\!\not|q-1$,}
\\
\left(\matrix 1&0\\a(\rho)\epsilon&1\endmatrix\right),&\text{otherwise.}\endcases$$
for some $a(\rho)\in A$.
\endproclaim
\definition{Proof} Let $\pi:GL_2(A[\epsilon])\rightarrow GL_2(A)$ denote the homomorphism induced by the augmentation map $a:A[\epsilon]\rightarrow A$. Because $\pi\circ\rho(i_{\goth p})$ is the identity matrix we have:
$$\rho(i_{\goth p})=\left(\matrix 1+a\epsilon&b\epsilon\\c\epsilon&
1+d\epsilon\endmatrix\right)$$
for some $a,b,c,d\in A$. The determinant of $\rho(i_{\goth p})$ is trivial hence $a+d=0$. When $l$ does not divide $q-1$ the matrix $\rho(i_{\goth p})$ fixes the vector $(1,1)$ hence $a+b=c+d=0$. When $l$ does divide $q-1$ the matrix $\rho(i_{\goth p})$ fixes the vector $(0,1)$ hence $b=d=0$. The claim is now clear.\ $\square$
\enddefinition
Let $l(\goth p)$ denote the largest power of $l$ dividing $N(\goth p)$.
\proclaim{Proposition 8.7} Assume that $\rho$ is in standard form. Using the notation of the lemma above we have $l(\goth p)a(\rho)=0$.
\endproclaim
\definition{Proof} For the sake of simple notation let $d$ denote $\deg(\goth p)$. Assume first that $l$ does not divide $q-1$. Then
$$\split\left(\matrix 1+q^da(\rho)\epsilon&-q^da(\rho)\epsilon
\\q^da(\rho)\epsilon&1-q^da(\rho)\epsilon\endmatrix\right)=&
\rho(i_{\goth p})^{q^d}=\rho(\text{Fr}_{\infty})^{-d}
\rho(i_{\goth p})\rho(\text{Fr}_{\infty})^d\\=&\left(\matrix 1+a(\rho)\epsilon&-q^da(\rho)\epsilon\\q^{-d}a(\rho)\epsilon&
1-a(\rho)\epsilon\endmatrix\right)\endsplit$$
by Lemma 8.4 hence $(q^d-1)a(\rho)=0$. On the other hand $\prod_{n=0}^{d-1}(\text{Fr}_{\infty}^{-n}i_{\goth p}\text{Fr}_{\infty}^n)^{q^{-n}}$ lies in the image of the inertia group $I_{\infty}$ in $\bold G_l$ by Lemma 8.4 hence the matrix:
$$\split\prod_{n=0}^{d-1}(\rho(\text{Fr}_{\infty})^{-n}
\rho(i_{\goth p})\rho(\text{Fr}_{\infty})^n)^{q^{-n}}=&\prod_{n=0}^{d-1}
\left(\matrix 1+a(\rho)\epsilon&-q^na(\rho)\epsilon\\q^{-n}a(\rho)\epsilon&
1-a(\rho)\epsilon\endmatrix\right)^{q^{-n}}\\=&
\left(\matrix 1+\sum_{n=0}^{d-1}q^{-n}
a(\rho)\epsilon&-da(\rho)\epsilon\\\sum_{n=0}^{d-1}q^{-2n}a(\rho)\epsilon&
1-\sum_{n=0}^{d-1}q^{-n}a(\rho)\epsilon\endmatrix\right)\endsplit$$
is upper-triangular with ones on the diagonal. Hence $(q^d-1)/(q-1)$ and $(q^{2d}-1)/(q^2-1)$ annihilates $a(\rho)$, too. The greatest common divisor of these numbers is:
$$\bigg({q^d-1\over q-1},{q^{2d}-1\over q^2-1}\bigg)=N(\goth p)\cdot
\cases(1,{q^d+1\over q+1}),&\text{if $d$ is odd,}
\\(q+1,q^d+1),&\text{otherwise.}\endcases$$
In any case the greatest common divisor on the right hand side of the equation above  divides $2$. Because $l$ does not divide $q-1$ it is odd so the claim above holds in this case. Assume now that $l$ does divide $q-1$. Then
$$\split\left(\matrix1&0\\q^da(\rho)\epsilon&1\endmatrix\right)=&
\rho(i_{\goth p})^{q^d}=\rho(\text{Fr}_{\infty})^{-d}
\rho(i_{\goth p})\rho(\text{Fr}_{\infty})^d\\=&\left(\matrix 1+{q^{-d}-1\over q-1}a(\sigma)\epsilon&{(q^{-d}-1)(q^d-1)\over(q-1)^2}
a(\sigma)\epsilon\\q^{-d}a(\rho)\epsilon&
1-{q^{-d}-1\over q-1}a(\sigma)\epsilon\endmatrix\right)\endsplit$$
hence $(q^d-1)/(q-1)$ annihilates $a(\rho)$. As we saw above the matrix:
$$\split\prod_{n=0}^{d-1}\rho(\text{Fr}_{\infty}^{-n}
i_{\goth p}\text{Fr}_{\infty}^n)^{q^{-n}}=&\prod_{n=0}^{d-1}
\left(\matrix 1+{q^{-n}-1\over q-1}a(\sigma)\epsilon&{(q^{-n}-1)(q^n-1)\over(q-1)^2}
a(\sigma)\epsilon\\q^{-n}a(\rho)\epsilon&
1-{q^{-n}-1\over q-1}a(\sigma)\epsilon\endmatrix\right)^{q^{-n}}\\=&
\left(\matrix1+\sum_{n=0}^{d-1}q^{-n}{q^{-n}-1\over q-1}a(\rho)\epsilon&-\sum_{n=0}^{d-1}\big({q^{-n}-1\over q-1}\big)^2a(\rho)\epsilon\\
\sum_{n=0}^{d-1}q^{-2n}a(\rho)\epsilon&1-\sum_{n=0}^{d-1}q^{-n}{q^{-n}-1\over q-1}a(\rho)\epsilon\endmatrix\right)
\endsplit$$
is upper-triangular with ones on the diagonal hence $(q^{2d}-1)/(q^2-1)$ annihilates $a(\rho)$ in this case, too. So the claim holds unless $d$ and $l$ are even. In the latter case we also need to use that the number:
$$q^{2d-2}\bigg(\sum_{n=0}^{d-1}q^{-n}{q^{-n}-1\over q-1}\bigg)=
{1\over q-1}
\bigg({q^{2d}-1\over q^2-1}-q^{d-1}{q^d-1\over q-1}\bigg)
=-{q^d-1\over q^2-1}\cdot{q^{d-1}-1\over q-1}$$
annihilates $a(\rho)$. Since
$${q^{d-1}-1\over q-1}=q^{d-2}+\cdots+1\equiv d-1\equiv1\mod2$$
in this case, the claim holds when $l=2$ and $d$ is even as well.\ $\square$
\enddefinition
\proclaim{Proposition 8.8} The group $I_R/I_R^2$ is cyclic and its order divides $l(\goth p)$.
\endproclaim
\definition{Proof} It will be sufficient to prove the same about the group $(I_R,l^n)/(I_R^2,l^n)$ for every positive integer $n$. Because $(I_R,l^n)/(I_R^2,l^n)$ is finite and annihilated by $l^n$ it will be sufficient to prove the same about the group $\text{Hom}((I_R,l^n)/(I_R^2,l^n),\Bbb Z/(l^n))$. Let $A$ be an object of $\Cal C$. We say that a $\Bbb Z_l$-homomorphism $\psi:R(\goth p)\rightarrow A[\epsilon]$ is $I_R$-admissible if $\psi(I_R)$ lies in the ideal of $A[\epsilon]$ generated by $\epsilon$. Note that a $\Bbb Z_l$-homomorphism $\psi:R(\goth p)\rightarrow A[\epsilon]$ is $I_R$-admissible if and only if the element of $\text{Def}_{\goth p}(A[\epsilon])$ corresponding to $\psi$ lies in $\text{Def}^0_{\goth p}(A[\epsilon])$.
For every $I_R$-admissible $\Bbb Z_l$-homomorphism $\psi:R(\goth p)\rightarrow
\Bbb Z/(l^n)[\epsilon]$ let $T(\psi):(I_R,l^n)\rightarrow\Bbb Z/(l^n)$ denote the unique $\Bbb Z_l$-linear map such that $\psi(i)=T(\psi)(i)\epsilon$ for every $i\in(I_R,l^n)$. The homomorphism $T(\psi)$ factors through $(I^2_R,l^n)$ and the map given by the rule $\psi\mapsto T(\psi)$ induces a bijection $b_n$ between the set $\text{Adm}(\Bbb Z/(l^n))$ of $I_R$-admissible $\Bbb Z_l$-homomorphisms $\psi:R(\goth p)\rightarrow\Bbb Z/(l^n)[\epsilon]$ and the set $\text{Hom}((I_R,l^n)/(I_R^2,l^n),\Bbb Z/(l^n))$. Under this bijection the set $\text{Adm}(\Bbb Z/(l^n))$ is equipped with a group structure.

Let $\rho_1$ and $\rho_2$ be two representations in standard form representing two elements of the set $\text{Def}^0_{\goth p}(\Bbb Z/(l^n)[\epsilon])$. We define the sum $\rho_3$ of $\rho_1$ and $\rho_2$ as the unique representation in standard form such that the latter, as a representation of $\text{Gal}(\overline F|F)$ on $\Bbb Z/(l^n)[\epsilon]^2$, is given by the rule:
$$\rho_3(g)=\rho_1(g)+\rho_2(g)-\rho_0(g),\quad\forall g\in\text{Gal}(\overline F|F),$$
where $\rho_0$ is the representation $\rho_{A[\epsilon]}$ constructed in the proof of Lemma 7.3. It is very easy to see that this operation makes the set of representations in standard form into a group $\text{St}(\Bbb Z/(l^n))$ with $\rho_0$ as the zero element. Let $\psi_i$ be the element of $\text{Adm}(\Bbb Z/(l^n))$ corresponding to $\rho_i$ where $i=1,2,3$. Then $\psi_3$ is the sum of $\psi_1$ and $\psi_2$. Hence $\text{Hom}((I_R,l^n)/(I_R^2,l^n),\Bbb Z/(l^n))$ is the quotient of $\text{St}(\Bbb Z/(l^n))$. Let $\rho$ be an element of $\text{St}
(\Bbb Z/(l^n))$ such that $a(\rho)=0$ using the notation of Lemma 8.6. Then $\rho$ is unramified hence it represents the unique element of $\text{Def}_{min}(\Bbb Z/(l^n)[\epsilon])$ which corresponds to the zero element of $\text{Adm}(\Bbb Z/(l^n))$. The claim now follows from Proposition 8.7.\ $\square$
\enddefinition
\proclaim{Corollary 8.9} The Eisenstein ideal $\goth E(\goth p)$ is locally principal.
\endproclaim
\definition{Proof} Note that the maximal ideal of the local ring $R(\goth p)$ is $(I_R,l)$. According to Propositions 8.8 the $\Bbb F_l$-dimension of the reduced tangent space $(I_R,l)/(I_R^2,l)$ is at most one hence $R(\goth p)$ is generated by a single element over $\Bbb Z_l$. Because $T(\goth p)$ is a quotient of $R(\goth p)$ by Proposition 7.8 the latter is also generated by a single element over $\Bbb Z_l$.\ $\square$
\enddefinition
\definition{Definition 8.10} Let us forget for a moment the notation introduced in chapter seven and let $\Cal O$ be a complete discrete valuation ring with residue field $\bold k$. Suppose that we have a commutative diagram of surjective homomorphism of complete Noetherian local $\Cal O$-algebras:
$$\diagram R& &\rTo^{\phi}& &T\\
&\rdTo_{\pi_R}& &\ldTo_{\pi_T}& \\
& &\Cal O.& &\\\enddiagram$$
Assume that $T$ is finitely generated and free as an $\Cal O$-module. Let $I_R$ and $I_T$ denote the kernel of the homomorphism $\pi_R$ and $\pi_T$, respectively. The congruence ideal $\eta_T$ of $T$ is defined to be the image $\pi_T(\text{Ann}_T(I_T))$ of the annulator of the ideal $I_T$ in the algebra $T$. For every finitely generated and torsion $\Cal O$-module $M$ let $\text{length}_{\Cal O}(M)$ denote the number of Jordan-H\"older components of $M$.
\enddefinition
\proclaim{Theorem (Wiles-Lenstra) 8.11} Suppose that $\eta_T\neq0$. Then
$$\text{\rm length}_{\Cal O}(I_R/I_R^2)\leq\text{\rm length}_{\Cal O}(\Cal O/\eta_T)$$
and equality holds if and only if $\phi$ is an isomorphism between complete intersections over $\Cal O$.
\endproclaim
\definition{Proof} This is Criterion I of [3] on page 343.\ $\square$
\enddefinition
\proclaim{Theorem 8.12} The $\Bbb Z_l$-homomorphism $\phi:R(\goth p)\rightarrow T(\goth p)$ is an isomorphism.
\endproclaim
Note that this result with Proposition 8.8 and Theorem 4.5 implies Theorem 1.2.
\definition{Proof} We wish to use Theorem 8.11 when $\Cal O=\Bbb Z_l$, $R=R(\goth p)$, $T=T(\goth p)$ and $\phi$ is the homomorphism constructed in the proof of Proposition 7.8. Let $\pi_R$ denote the homomorphism defined in Definition 7.4: in this case $I_R$ is the ideal introduced in Definition 7.4. Let $I_T$ denote the ideal generated by the image of $\widehat{\goth E}(\goth p)$ with respect to the projection $\widehat{\Bbb T}(\goth p)\rightarrow T(\goth p)$. As we saw in the proof of Proposition 7.8 we have $\Bbb Z_l=T(\goth p)/I_T$ and the factorization map $\pi_T:T(\goth p)\rightarrow T(\goth p)/I_T$ makes the diagram in Definition 8.10 commutative.  Let $T_0(\goth p)$ denote the completion of $\Bbb T(\goth p)$ with respect to the maximal ideal $(l,\goth E(\goth p))$ and let $\pi_0:T(\goth p)\rightarrow T_0(\goth p)$ denote the surjection induced by the $\widehat{\Bbb T}(\goth p)\rightarrow\Bbb T(\goth p)$. The direct sum $\pi_0\oplus\pi_T:T(\goth p)\rightarrow T_0(\goth p)\oplus\Bbb Z_l$ is injective and $\text{Ann}_{T(\goth p)}(I_T)$ is equal to the part of $(\pi_0\oplus\pi_T)(T(\goth p))$ lying in the second factor $\Bbb Z_l$. If $0\oplus n\in\text{Ann}_{T(\goth p)}(I_T)$ then $n\oplus0\in(\pi_0\oplus\pi_T)(T(\goth p))$. The part of $\pi_0\oplus\pi_T(T(\goth p))$ lying in the first factor $T_0(\goth p)$ is the ideal generated by the image of $\goth E(\goth p)$ with respect to the canonical map $\Bbb T(\goth p)\rightarrow T_0(\goth p)$. According to claim $(vi)$ of Proposition 7.11 of [19] on pages 180-181 the cyclic group $\Bbb T_l(\goth p)/\goth E_l(\goth p)$ is finite and its order is divisible by $l(\goth p)$. Hence $n\Bbb Z_l\subseteq l(\goth q)\Bbb Z_l$ for every $n$ as above and the criterion of Wiles-Lenstra is satisfied by Proposition 8.8.\ $\square$
\enddefinition

\heading 9. Diophantine applications
\endheading

\definition{Definition 9.1} The deformation theoretical methods of this paper also give an alternative route to prove the main diophantine results of the paper [19]. The most important advantage of this method is that we do not need the analysis of the special fiber of the modular curves $X_0(\goth p)$ and $X_1(\goth p)$ at the prime $\goth p$ which occupies chapter 8 of the paper quoted above. In the rest of this paper we give a short description of this alternative route. We will repeat some of the arguments presented in [19] for the sake of exposition. Let $\Cal E(\goth p)$ denote the largest torsion subgroup of $J_0(\goth p)(\overline F)$ annihilated by the Eisenstein ideal $\goth E(\goth p)\triangleleft\Bbb T(\goth p)$. In the proof of the following key proposition the Gorenstein property plays an essential role.
\enddefinition
\proclaim{Proposition 9.2} The order of the group $\Cal E(\goth p)$ divides $N(\goth p)^2$.
\endproclaim
\definition{Proof} For every group $G$ and prime number $l$ let $G_l$ denote its maximal $l$-primary subgroup. It is sufficient to prove that $\Cal E(\goth p)_l$ divides $l(\goth p)^2$ for every prime number $l$ where $l(\goth p)$ denote the largest power of $l$ dividing $N(\goth p)$. By Theorem 1.4 we may assume that $l\neq p$. Let $\mu_{\infty}\subset\Bbb C^*_{\infty}$ denote the subgroup of the roots of unity as in the proof of Lemma 5.9. The Pontryagin dual $P_l$ of the $l$-primary torsion $Q_l$ of the group Hom$(\overline{\Gamma}_0(\goth p),\mu_{\infty})$ is the $\Bbb Z_l$-dual of a locally free $\Bbb T_l(\goth p)$-module of rank one by claim $(v)$ of Proposition 7.11 of [19] on pages 160-161. In particular the part of $Q_l$ annihilated by the Eisenstein ideal has order $l(\goth p)$ by the Gorenstein property. On the other hand the Pontryagin dual of the $l$-primary torsion $R_l$ of the quotient of the torsion group of $J_0(\goth p)(\Bbb C_{\infty})$ by the image of Hom$(\overline{\Gamma}_0(\goth p),\mu_{\infty})$ with respect to the map $\Psi_{AJ}$ is the $\Bbb Z_l$-dual of $P_l$. Hence the part of $R_l$ annihilated by the Eisenstein ideal has order $\goth l(\goth p)$, too. Therefore $\Cal E(\goth p)_l$ has a two-step filtration such that the orders of the steps divide $l(\goth p)$.\ $\square$
\enddefinition
Let $\Cal C(\goth p)\subset J_0(\goth p)(\overline F)$ and $\Cal S(\goth p)\subset
J_0(\goth p)(\overline F)$ denote the cuspidal divisor group and the Shimura group, respectively. (For definition see Definition 1.3 and Definition 1.5 of [19] on page 132, respectively.) The group $\Cal C(\goth p)$ lies in $J_0(\goth p)(F)$ while the group $\Cal S(\goth p)$ is $\mu$-type. Let $t(\goth p)$ denote the greatest common divisor of $N(\goth p)$ and $q-1$. We will need the following facts:
\proclaim{Proposition 9.3} The following holds:
\roster
\item"$(i)$" the groups $\Cal C(\goth p)$ and $\Cal S(\goth p)$ are cyclic of order $N(\goth p)$,
\item"$(ii)$" the intersection of $\Cal C(\goth p)$ and $\Cal S(\goth
p)$ is their unique cyclic group of order $t(\goth p)$.
\endroster
\endproclaim
\definition{Proof} Claim $(i)$ is just Corollary 5.11 of [7] on page 235 and Lemma 8.17 of [19] on page 171. For the proof of claim $(ii)$ see Proposition 9.3 of [19] on pages 175-176.\ $\square$
\enddefinition
In the paper [19] we also introduced another group $\Cal D(\goth p)[l]
\subset J_0(\goth p)(\overline F)$ for every prime number $l$ dividing $t(\goth p)$. (For definition see  Definitions 9.4 and 9.17 on pages 176 and 183 of [19], respectively.) The group $\Cal D(\goth p)[l]$ is $l$-torsion and contains the $l$-torsion subgroup $\Cal S(\goth p)[l]$ of the Shimura group $\Cal S(\goth p)$. Let $\Cal F(\goth p)[l]$ denote the quotient of $\Cal D(\goth p)[l]$ by its subgroup $\Cal S(\goth p)[l]$.
\proclaim{Proposition 9.4} The following holds:
\roster
\item"$(i)$" the Galois modules $\Cal S(\goth p)[l]$ and $\Cal F(\goth
p)[l]$ are constant of order $l$,
\item"$(ii)$" the Galois module $\Cal D(\goth p)[l]$ is everywhere
unramified,
\item"$(iii)$" the Galois module $\Cal D(\goth p)[l]$ is contained in
$\Cal E(\goth p)$,
\item"$(iv)$" the exact sequence:
$$0@>>>\Cal S(\goth p)[l]@>>>\Cal D(\goth p)[l]@>>>\Cal F(\goth
p)[l]@>>>0$$
of Galois modules does not split over $F$.
\endroster
\endproclaim
\definition{Proof} This is just Proposition 9.18 of [19] on pages 184-185.\ $\square$
\enddefinition
Let $\Cal T(\goth p)\subset J_0(\goth p)(\overline F)$ and $\Cal M(\goth p)\subset
J_0(\goth p)(\overline F)$ denote the torsion subgroup of $J_0(\goth p)(F)$ and the maximal $\mu$-type \'etale subgroup scheme of $J_0(\goth p)$, respectively. The following lemma is an easy consequence of the fact that neither $\Cal T(\goth p)$ nor $\Cal M(\goth p)$ has $p$-torsion and the Eichler-Shimura relations.
\proclaim{Lemma 9.5} The groups $\Cal T(\goth p)$ and $\Cal M(\goth p)$ are contained by $\Cal E(\goth p)$.
\endproclaim
\definition{Proof} For proof see Lemma 7.16 on page 163 and Lemma 10.4 on pages 186-187 of the paper [19].\ $\square$
\enddefinition
The main diophantine result of [19] is:
\proclaim{Theorem 9.6} We have $\Cal T(\goth p)=\Cal C(\goth p)$ and $\Cal M(\goth p)=\Cal S(\goth p)$.
\endproclaim
\definition{Proof} We are going to show the equalities $\Cal T(\goth p)_l=\Cal C(\goth p)_l$ and $\Cal M(\goth p)_l=\Cal S(\goth p)_l$ for every prime number $l$. We may assume that $l$ is Eisenstein that is $l$ divides $N(\goth p)$. First assume that $l$ is a prime number which does not divide $q-1$. In this case the intersection of the
groups $\Cal C(\goth p)_l$ and $\Cal S(\goth p)_l$ is trivial because the action of the absolute Galois group of $F$ does not fix any non-zero element of an $l$-primary $\mu$-type Galois module. On the other hand the order of these groups is $l(\goth p)$ by claim $(i)$ of Proposition 9.3 hence their direct sum has the same order as $\Cal E(\goth p)_l$ by Proposition 9.2. Therefore $\Cal E(\goth p)_l$ and the direct sum of $\Cal C(\goth p)_l$ and $\Cal S(\goth p)_l$ are equal. Hence $\Cal C(\goth p)_l$ and $\Cal S(\goth p)_l$ is the maximal constant and $\mu$-type Galois module of $\Cal E(\goth p)_l$, respectively. So $\Cal T(\goth p)_l=\Cal C(\goth p)_l$ and $\Cal M(\goth p)_l=\Cal S(\goth p)_l$ by Lemma 9.5.

Assume now that $l$ is a prime number which does divide $q-1$. Then it also divides $t(\goth p)$ because it is Eisenstein. By claim $(iii)$ of Proposition 9.4 the group $\Cal E(\goth p)[l]$ contains $\Cal D(\goth p)[l]$. We saw in the proof of Proposition 9.2 that $\Cal E(\goth p)_l$ has a two-step filtration such that the steps are cyclic groups. Hence the orders of $\Cal E(\goth p)[l]$ and $\Cal D(\goth p)[l]$ are equal by claim $(i)$ of Proposition 9.4 so the groups $\Cal E(\goth p)[l]$ and $\Cal D(\goth p)[l]$ must be equal, too. Now we are going to prove that $\Cal T(\goth p)_l=\Cal C(\goth p)_l$ and $\Cal M(\goth p)_l=\Cal S(\goth p)_l$. If this claim is false then there is an element $x$ in $\Cal M(\goth p)_l-\Cal S(\goth p)_l$ (resp$\text{.}$ in $\Cal T(\goth p)_l-\Cal C(\goth p)_l$) such that $lx$ is in $\Cal S(\goth p)_l$ (resp$\text{.}$ in $\Cal C(\goth p)_l$). The element $x$ is annihilated by $l(\goth p)$, since it is annihilated by the Eisenstein ideal. Therefore $lx$ is annihilated by ${l(\goth p)\over l}$. Since both $\Cal S(\goth p)_l$ and $\Cal C(\goth p)_l$ are cyclic of order $l(\goth p)$, the element $lx$ must have an $l$-th root $u$ in $\Cal S(\goth p)_l$ (resp$\text{.}$ in
$\Cal C(\goth p)_l$) by the above. Subtracting $u$ from $x$ we get that
we may assume that $x$ is $l$-torsion. Hence we must have
$x\in\Cal D(\goth p)[l]$ by the above. Since $\Cal S(\goth p)[l]$ is the maximal constant as well as $\mu$-type subgroup of $\Cal D(\goth p)[l]$ by claims $(ii)$ and $(iv)$ of Proposition 9.4 we conclude that $x$ is actually in $\Cal S(\goth p)[l]$. The intersection of $\Cal S(\goth
p)_l$ and $\Cal C(\goth p)_l$ is exactly the largest constant Galois
submodule of the former by claim $(ii)$ of Proposition 9.3, so the claim is now
clear.\ $\square$
\enddefinition

\enddocument